\documentclass[a4paper,11pt]{article}
\usepackage[margin=1.2in]{geometry}
\pdfoutput=1
\usepackage{mathtools}

\usepackage[utf8]{inputenc}
\usepackage{amsfonts,amssymb, bm,amsthm}
\usepackage{graphicx,color,epstopdf}
\usepackage[inline]{showlabels}
\allowdisplaybreaks
\newtheorem{mylemma}{Lemma}[section]
\newtheorem{mytheorem}{Theorem}[section]

\newtheorem{myremark}{Remark}[section]
\newtheorem{mydef}{Definition}[section]

\def\XXint#1#2#3{{\setbox0=\hbox{$#1{#2#3}{\int}$}
    \vcenter{\hbox{$#2#3$}}\kern-.5\wd0}}

\def\tx {\tilde{x}}
\def\ty {\tilde{y}}
\def\tu {\tilde{u}}

\def\Log {{\rm Log}}

\def\rsc{{\rm sc}}

\def\wtd#1{{\widetilde{#1}}}

\def\cS {{\cal S}}
\def\cK {{\cal K}}
\def\cKd {{\cal K}'}
\def\cT {{\cal T}}
\def\hG {\hat{G}}
\def\hf {\hat{f}}

\def\bi{{\bf i}}

\begin{document}
\title{Mathematical analysis of wave radiation by a step-like surface} \author{
  Wangtao Lu$^1$ }
\footnotetext[1]{School of Mathematical Sciences, Zhejiang University, Hangzhou
  310027, China. Email: wangtaolu@zju.edu.cn}
\maketitle
\begin{abstract}
  This paper proposes, for wave propagating in a globally perturbed half plane
  with a perfectly conducting step-like surface, a sharp Sommerfeld radiation
  condition (SRC) for the first time, an analytic formula of the far-field
  pattern, and a high-accuracy numerical solver. We adopt the Wiener-Hopf method
  to compute the Green function for a cracked half plane, a background for the
  perturbed half plane. We rigorously show that the Green function
  asymptotically satisfies a universal-direction SRC (uSRC) and radiates purely
  outgoing at infinity. This helps to propose an implicit transparent boundary
  condition for the scattered wave, by either a cylindrical incident wave due to
  a line source or a plane incident wave. Then, a well-posedness theory is
  established via an associated variational formulation. The theory reveals that
  the scattered wave, post-subtracting a known wave field, satisfies the same
  uSRC so that its far-field pattern is accessible theoretically. For a
  plane-wave incidence, asymptotic analysis shows that merely subtracting
  reflected plane waves, due to non-uniform heights of the step-like surface at
  infinity, from the scattered wave in respective regions produces a
  discontinuous wave satisfying the uSRC as well. Numerically, we adopt a
  previously developed perfectly-matched-layer (PML) boundary-integral-equation
  method to solve the problem. Numerical results demonstrate that the PML
  truncation error decays exponentially fast as thickness or absorbing power of
  the PML increases, of which the convergence relies heavily on the 
  Green function exponentially decaying in the PML.
\end{abstract}
\section{Introduction}
\begin{figure}[!ht]
  \centering
  (a)\includegraphics[width=0.4\textwidth]{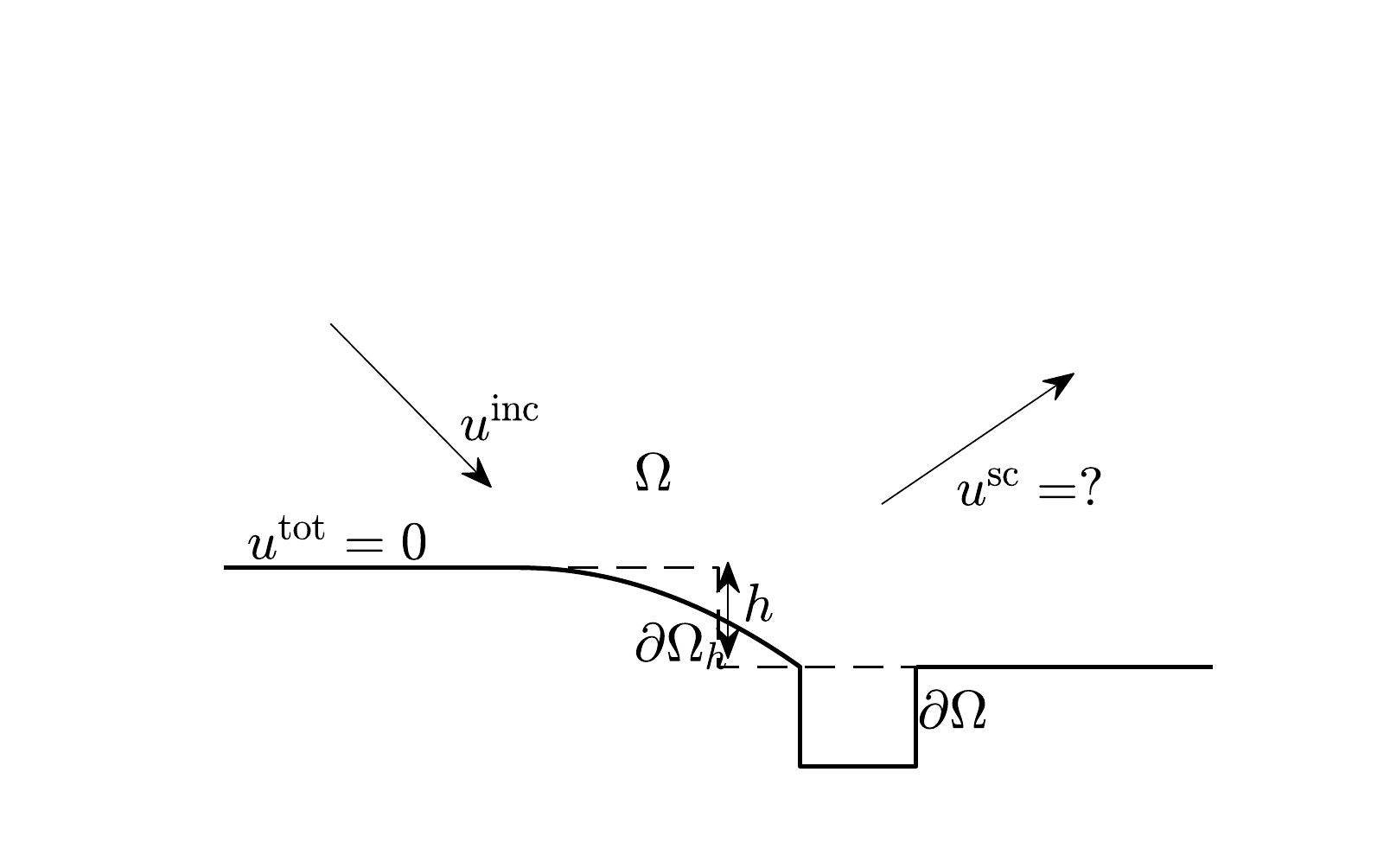}
  (b)\includegraphics[width=0.4\textwidth]{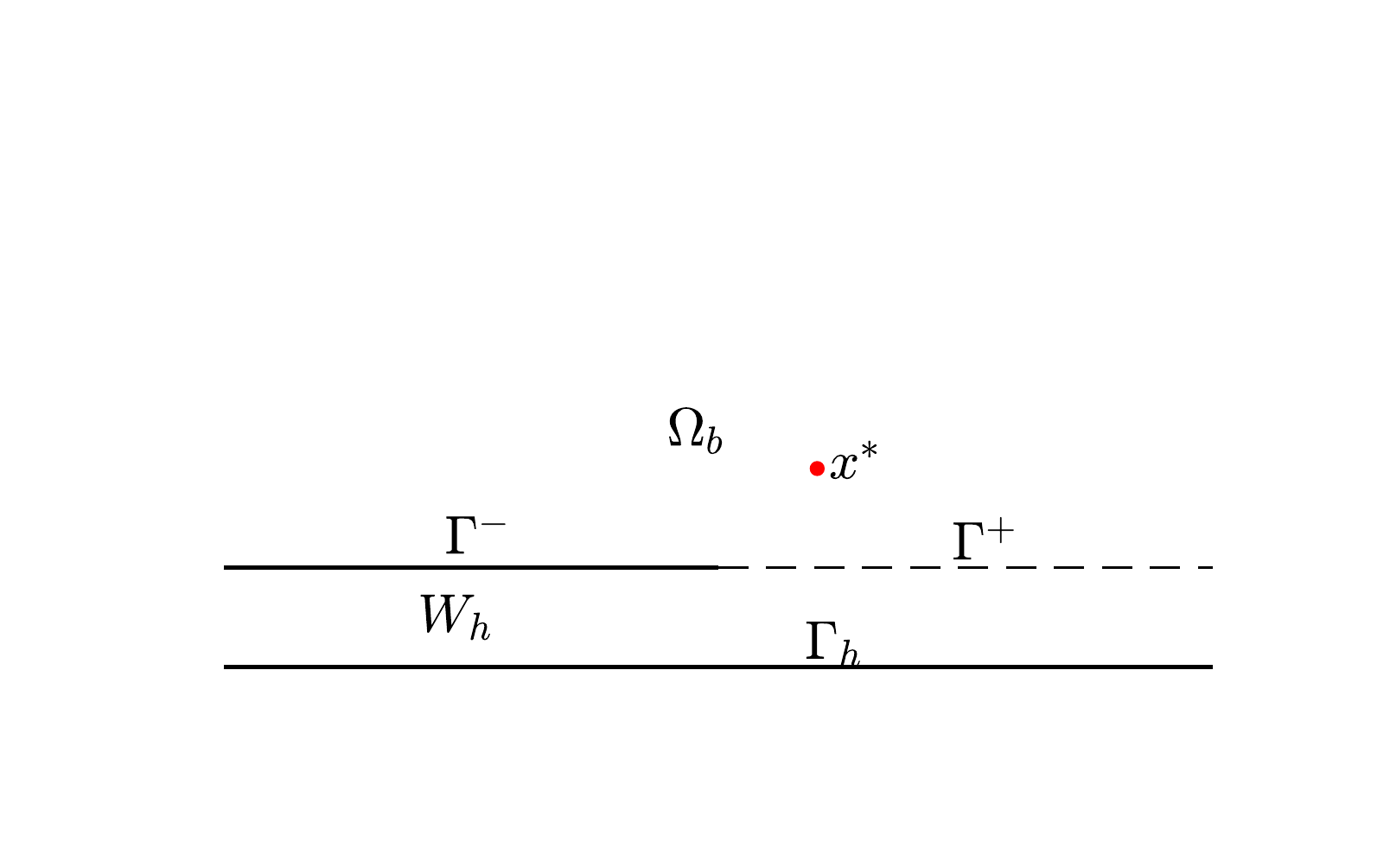}
  \caption{(a): A schematic of the scattering problem; the scattering
    surface $\partial\Omega$ is a local perturbation of the surface
    $\partial\Omega_h$ indicated by dashed lines. (b): A reference model, where 
  the dot in red indicates the location of a point source $x^*$; $\Omega_b$ can
  be regarded as a background domain for the scattering problem in (a). All
  solid lines representing perfectly conducting surfaces such that $u^{\rm tot}=0$ there.}
  \label{fig:model}
\end{figure}

Wave propagating in inhomogeneous media has numerous applications in both
scientific and engineering areas \cite{chew95,colkre13}. This paper
mathematically analyzes wave scattering in a globally perturbed two-dimensional (2D) half-space with a
perfectly conducting step-like surface. Specifically, the scattering problem is
setup as follows. Let an incident wave $u^{\rm inc}$ be specified in an
unbounded domain $\Omega\subset \mathbb{R}^2$. As illustrated in
Figure~\ref{fig:model}(a),
the scattering surface $\partial\Omega$, i.e., the boundary of $\Omega$, is
piecewise Lipschitz and is from locally perturbing a step-function curve
\[
  \partial\Omega_h=\{(x_1,0):x_1\leq 0\}\cup\{(0,x_2):-h<x_2<0\}\cup\{(x_1,-h):x_1\geq 0\}, 
\]
the boundary of $\Omega_h$, where $h$ denotes the height of the step. Moreover,
we suppose $\partial\Omega$ satisfies the following geometrical condition (c.f.
\cite{chaels10})
\[
  (x_1,x_2)\in\Omega\Rightarrow (x_1,x_2+a)\in\Omega,\quad\forall a\geq 0.
\]
The total wavefield $u^{\rm tot}$, the sum of incident wave $u^{\rm inc}$ and
the resulting scattered wave $u^{\rm sc}$, solves
\begin{align}
  \label{eq:gov:tot1}
  \Delta u^{\rm tot} + k^2u^{\rm tot} &= 0,\quad{\rm on}\quad\Omega,\\
  \label{eq:gov:tot2}
  u^{\rm tot} &= 0,\quad{\rm on}\quad\partial\Omega,
\end{align}
where the Laplace operator $\Delta=\partial_{x_1}^2 + \partial_{x_2}^2$,
$k=\frac{2\pi}{\lambda}$ denotes the wavenumber and $\lambda$ is the wavelength,
and zero Dirichlet boundary condition is imposed on boundary $\partial\Omega$.
Physically, $u^{\rm tot}$ could represent the longitudinal component of an
electric field due to a perfectly conducting (PEC) surface $\partial\Omega$, or
a sound wave due to a sound-soft surface $\partial\Omega$. Mathematically, a
fundamental question is how to prove the well-posedness of the problem.

To answer this question, a ``proper'' radiation condition of $u^{\rm tot}$ at
infinity should be posed first. On one hand, the proposed radiation condition
must be physically reliable for characterizing wave propagation in the medium;
on the other hand, it should be helpful in designing an accurate and efficient
numerical solver. In the past decades, existing literature have proposed, in
general, three types of radiation conditions responsible for the well-posedness
of wave scattering by such a rough surface.

The first one is the angular spectrum representation (ASR) condition
\cite{des02,baohuyin18}, a.k.a the upward propagation radiation condition (UPRC)
\cite{zhacha98,zhacha03}. This UPRC condition requires that the scattered wave
$u^\rsc$ contain no downgoing waves in its Fourier-integral representation by
plane waves, and has been successfully used for establishing the well-posedness
of scattering problems for more general rough surfaces
\cite{chaels10,chamon05,chazha99}, certainly including the scattering surface
$\partial\Omega$. The second one is the modal radiation condition (MRC)
\cite{bongouhaz11,bontil01a,bontil01b}, which is quite suitable for analyzing
wave scattering by rectangular interfaces at infinity, e.g. $\partial\Omega_h$.
The propagation domain can be decomposed into a finite number of rectangular
regions, and waves satisfying MRC can be expanded as the sum of outgoing
eigenmodes in any exterior rectangular region; existing numerical mode matching
methods \cite{biederbaeolydez01, lushilu14, luluson19} are highly related to the
MRC condition. Probably, the most attractive one is the well-known Sommerfeld
radiation condition (SRC) pioneered by Sommerfeld \cite{som49}. Roughly
speaking, a wave satisfying SRC propagates purely outgoing at infinity so that
SRC is the strongest condition among the aforementioned three conditions. Thanks
to this important feature, SRC and its derivatives have been considerably used
for well posing various wave propagation problems \cite{li10jcp,
  roazha92,baohuyin18}. Its advantage is threefold. Firstly, by use of a known
background Green function, SRC helps to give an exterior Green's representation
formula, which in turn determines an exact transparent boundary condition (TBC)
to terminate the unbounded domain. Secondly, the exterior Green's representation
formula together with the far-field pattern of the background Green function,
directly determines the far-field pattern of any outgoing wave, which plays an
important role in both direct and inverse problems \cite{colkre13}. Thirdly,
from a numerical point of view, the outgoing behavior makes absorbing boundary
conditions, such as perfectly matched layer (PML) \cite{ber94,chewee94},
applicable to numerically truncate the unbounded domain such that standard
numerical methods, such as finite element method \cite{mon03} or boundary
integral equation (BIE) method \cite{luluqia18,hulurat20}, could then apply.
Bonnet-Bendhia {\it et al.} \cite{bonfliton18} have recently proposed a
half-space matching method to solve scattering problems in infinite media,
discretizing the Fourier-integral representation of an outgoing wave along lines
enclosing scatters, which have close relations with the UPRC and SRC conditions.

For a cylindrical-wave incidence due to a line source, our previous work
 \cite{hulurat20} has shown that for a general rough surface, the scattered wave $u^\rsc$ directly satisfies the integral
form of SRC, a weakened SRC (wSRC); we also show that if a horizontal stripe is
removed from the scattering domain, then $u^\rsc$ satisfies the classic SRC, a
pointwise but stronger condition,   in
the remaining half plane. But for a plane-wave incidence, no SRC condition has
been imposed for wave scattered by the unbounded curve $\partial\Omega$, until
recently the author and Hu in \cite{luhu19} proved that the scattered wave
$u^\rsc$ piecewisely satisfies wSRC at infinity. Due to the absence of a
background Green function, the wSRC condition couldn't strictly characterize the
radiation behavior of $u^{\rm sc}$ at infinity, and hence the far field pattern
of $u^\rsc$ is unclear; moreover, we couldn't explain why the PML truncation
used there should produce a physically correct solution. Motivated by this, this
paper makes use of the Green function of a special background domain, to
rigorously derive for the scattered wave $u^\rsc$, a sharper SRC condition, a
closed form of its far field pattern, and a high-accuracy numerical solver.

Instead of regarding $\Omega_h$ in Figure~\ref{fig:model}(a) as the background,
we choose the cracked half-plane $\Omega_b$, as illustrated in
Figure~\ref{fig:model}(b). As inspired by \cite{mitlee71, tan99}, the Green
function $G$ for such a background could be analytically computed based on the
well-known Wiener-Hopf technique \cite{nob58,gak66}. We point out that
$\Omega_h$ is not a suitable background since the vertical segment of
$\partial\Omega_h$ makes Fourier transforming $x_1$-variable impossible. We
rigorously analyze asymptotic behavior and far-field pattern of the Green
function $G$ at infinity, showing that it satisfies a universal-direction SRC
(uSRC) stronger than that in \cite{luhu19,hulurat20} in the sense that $G$
satisfies the classic SRC condition uniformly along all directions in the full
region $\Omega_b$, with no restrictions on the satisfied region. Based on this,
we propose, for either a cylindrical incident wave due to a line source or a
plane incident wave, an implicit TBC (ITBC) to terminate the unbounded domain
$\Omega$ such that a well-posedness theory is established via an associated
variational formulation. The theory reveals that the scattered wave,
post-subtracting a known wave field, satisfies the same uSRC, so that its far-field
pattern is accessible theoretically. For a plane-wave incidence, asymptotic
analysis shows that merely subtracting reflected plane waves, due to non-uniform
heights of the step-like surface at infinity, from the scattered wave in
respectively affected regions, as was done in the previous paper \cite{luhu19},
produces a radiating wave, though discontinuous, satisfies the uSRC as well.
Numerically, we adopt a previously developed PML-based BIE method
\cite{luluqia18} to solve the problem. Numerical results demonstrate that the
truncation error due to PML decays exponentially fast as thickness or absorbing
power of the PML increases, of which the convergence relies heavily on the
outgoing Green function $G$ decaying exponentially in the PML.

The rest of this paper is organized as follows. In section 2, we derive present
a closed-form of the Green function for the background $\Omega_b$ derived by the
Wiener-Hopf technique and analyze its asymptotic behavior at infinity. In
section 3, we propose for $u^\rsc$ the uSRC condition, introduce the ITBC
condition, pose for the scattering problem a variational formulation, and
establish the well-posedness theory. We present the PML-BIE method to
numerically solve the scattering problem and briefly analyze the property of
Green function in the PML in section 4, and draw our conclusion finally.
\subsection{Some notations}
We shall adopt the notations of Sobolev spaces described in \cite{mcl00}. For a
generic bounded Lipschitz domain $B$ with boundary $\partial B$, let $H^1(B)$
and $H^{1/2}(\partial B)$ be standard domain and boundary Sobolev spaces
equipped with norms $||\cdot||_{H^1(B)}$ and $||\cdot||_{H^{1/2}(\partial B)}$,
and $\widetilde{H^{-1}}(B)=H^1(B)^*$ and $H^{-1/2}(\partial B)=H^{1/2}(\partial
B)^*$ be the associated dual spaces, respectively. We shall use
$(\cdot,\cdot)_{B}$ to denote the standard inner product in $L^2(B)$, and
$<\cdot,\cdot>_B$ to denote the duality pairing between two dual spaces on $B$,
etc.. Let $C_{\rm comp}^{\infty}(B)$ be the space of smooth functions with
compact support in $B$. Let $\gamma: H^1(B)\to H^{1/2}(\partial B)$ be the
bounded trace operator. We introduce Sobolev spaces on partial boundaries of
$B$. Let $\partial B_1\subsetneq \partial B$ be Lipschitz, then
\begin{align*}
  H^{1/2}(\partial B_1) &= \{u: u=U|_{\partial B}, U\in H^{1/2}(\partial B)\},\\
  \wtd{H^{1/2}}(\partial B_1) &= \{U\in H^{1/2}(\partial B): {\rm supp}U\ \subset \overline{\partial B_1}\}.
\end{align*}
Thus, $\wtd{H^{-1/2}}(\partial B_1)=H^{1/2}(\partial B_1)^*$ and
$H^{-1/2}(\partial B_1)=\wtd{H^{1/2}}(\partial B_1)^*$ are the dual spaces. We shall
need the following bounded extension operator $E: H^{-1/2}(\partial B_1)\to
H^{-1/2}(\partial B)$ such that $E\phi|_{\partial B_1}=\phi$ for any
$\phi\in H^{-1/2}(\partial B_1)$ and then $E^*$ denotes the adjoint bounded operator $E^*:
{H^{1/2}}(\partial B)\to \wtd{H^{1/2}}(\partial B_1)$. It can be seen that $E^*$
restricted on $\wtd{H^{1/2}}(\partial B)$ becomes identical.

For $R>0$, let $D_R$ denotes the disk of radius $R$ centered at origin. For the
unbounded domain $\Omega$, let $H^{1}_{\rm loc}(\Omega)$ be
the space of elements in $H^{1}(\Omega \cap D_R)$ for any $R>0$.

\section{The Green function of a cracked half-plane}
We consider the following auxiliary problem
\begin{align}
  \label{eq:model}
  \Delta G(x;x^*) + k^2 G(x;x^*) &= -\delta(x-x^*),\quad {\rm on}\quad \Omega_b,\\
  \label{eq:model:bc}
  G &= 0,\quad{\rm on}\quad \Gamma_h\cup\Gamma^-,
\end{align}
where $\Omega_b=\mathbb{R}_2^+\cup W_h\cup \Gamma^+$,
$W_h=\{(x_1,x_2):-h<x_2<0\}$, $\Gamma^+=\{(x_1,0):x_1>0\}$, a PEC surface
$\Gamma_h=\{(x_1,-h):x_1\in\mathbb{R}\}$ and a PEC crack
$\Gamma^-=\{(x_1,0):x_1<0\}$, $x=(x_1,x_2)$ and the source point
$x^*=(x_1^*,x_2^*)$, as illustrated in Figure~\ref{fig:model}(b). By the
Wiener-Hopf method, \cite{mitlee71, tan99} have computed the Green function $G$
excited by source $x^*$ in the waveguide $W_h$. Following the same
approach closely, we first give the closed-form of $G$ for any $x^*\in\Omega_b$
in this section, and shall briefly present the Wiener-Hopf method in Appendix
for the sake of completeness. Next, we rigorously analyze the asymptotic
behavior of $G$ at infinity, which, as we shall see, accounts for the asymptotic
behavior of $u^{\rm tot}$ of problem (\ref{eq:gov:tot1}) and
(\ref{eq:gov:tot2}).

\subsection{Closed-form of $G$}

We distinguish three cases: $x^*\in\mathbb{R}_+^2$, $\Gamma^+$ or $W_h$. For $x^*\in\mathbb{R}_+^2$ with $x_2^*>0$, the Green
function $G$ takes the following form
\begin{equation}
  \label{eq:G:1}
  G(x;x^*) =\left\{
    \begin{array}{ll}
      G_1(x;x^*)+G^{\rm in}(x;x^*), & x \in \mathbb{R}_2^+,\\
      G_2(x;x^*), & x \in W_h,
    \end{array}
  \right. 
\end{equation}
In the above, the incident wave $G^{\rm in}(x;x^*)= \Phi^k(x;x^*) - \Phi^k(x;x_{\rm im}^*)$,
where $x_{\rm im}^*=(x_1^*,-x_2^*)$ is the image point of $x^*$ about $x_2=0$, $\Phi^k$
denotes the free-space Green's function of wavenumber $k$, i.e.,
\[
  \Phi^k(x;x^*) := \frac{\bi}{4}H_0^{(1)}(k|x-x^*|) =
  \frac{\bi}{4\pi}\int_{-\infty}^{+\infty}\frac{e^{-\bi \xi (x-x_1^*) + \bi \mu |x_2-x_2^*|}}{\mu} d\xi,
\]
and $\mu=\sqrt{k^2-\xi^2}$ and throughout this paper, the branch cut of
$\sqrt{\cdot}$ is chosen as the negative real axis to make its real part
non-negative. The two scattered waves $G_1$ and $G_2$ are
\begin{align}
  \label{eq:sol:G1}
  G_1(x;x^*) =& \frac{1}{2\pi}\int_{\cal L} \hf^+(\xi;x^*)e^{-\bi \xi x_1+\bi \mu x_2}d\xi,\\
  \label{eq:sol:G2}
  G_2(x;x^*) =& \frac{1}{2\pi}\int_{\cal L}\hf^+(\xi;x^*) \frac{e^{-\bi\mu x_2}-e^{\bi\mu (x_2+2h) }}{1-e^{2\bi\mu h}} e^{-\bi \xi x_1}d\xi. 
\end{align}
Here, ${\cal L}$ denotes a smooth path from $-\infty+0\bi$ to $+\infty-0\bi$ in
the complex plane of $\xi$, passing through the origin $O$, slightly above the
negative real axis in quadrant $\mathbb{C}^{-+}$ and slightly below the positive
real axis in quadrant $\mathbb{C}^{+-}$ such that $1-e^{2\bi \mu h}$ is strictly
nonzero on ${\cal L}$; if $e^{2 \bi k h}=1$, then we could redefine ${\cal L}$
by slightly moving ${\cal L}$ leftward. For any $\xi\in{\cal L}$,
\begin{align}
  \label{eq:+t-}
  \hf^+(\xi;x^*)=&\frac{H^+(\xi;x^*)K^+(\xi)}{2\bi\sqrt{k+\xi}}=\frac{e^{\bi \mu x_2^*}(1-e^{2\bi\mu h})}{2\bi\mu}e^{\bi\xi x_1^*}-\frac{H^-(\xi;x^*)(1-e^{2\bi\mu h})}{2\bi\sqrt{k+\xi}K^-(\xi)}\\
  \label{eq:K+-:onL}
  K^{\pm}(\xi) =& \sqrt{1-e^{2\bi\mu(\xi) h}}\exp\left\{\pm \frac{1}{2\pi\bi}{\rm p.v.}\int_{\cal L} \frac{\Log(1-e^{2\bi\mu(t) h})}{t-\xi}dt\right\},\\
  \label{eq:H+-:onL}
  H^{\pm}(\xi;x^*) =& \frac{e^{\bi \xi x_1^*+\bi\mu x_2^*}K^-(\xi)}{2\sqrt{k-\xi}}\pm\frac{1}{2\pi\bi}{\rm p.v.}\int_{\cal L} \frac{e^{\bi tx_1^*+ \bi \mu(t) x_2^*}K^-(t)}{\sqrt{k-t}(t-\xi)}dt,
\end{align}
where p.v. indicates that the integral is a principal value integral, and the second
equality in (\ref{eq:+t-}) is based on the following two decompositions on
${\cal L}$,
\begin{align}
  \label{eq:decomp:1}
  1 - e^{2\bi \mu h} &= K^+(\xi) K^-(\xi),\\
  \label{eq:decomp:2}
  \frac{e^{\bi\xi x_1^*+\bi\mu x_2^*}}{\sqrt{k-\xi}}K^{-}(\xi) &= H^+(\xi;x^*) + H^-(\xi;x^*).
\end{align}

For $x^*\in\Gamma^+$ such that $x_2^*=0$ and $x_1^*>0$, we regard
$G(x;(x_1^*,0))$ as the limit of $G(x;(x_1^*,x_2^*))$ as $x_2^*\to 0^+$. In
doing so, $G$ remains invariant except that $H^{\pm}$ should be redefined as
follows,
\begin{equation}
    \label{eq:H0pm:L}
    H^{\pm}(\xi;(x_1^*,0))= \frac{1}{2}\frac{e^{\bi \xi x_1^*}(K^-(\xi)-1)}{\sqrt{k-\xi}}\pm \frac{1}{2\pi\bi}{\rm p.v.}\int_{\cal L}\frac{e^{\bi t x_1^*}(K^-(t)-1)}{\sqrt{k-t}(t-\xi)}dt +\left\{
      \begin{array}{l}
        \frac{e^{\bi \xi x_1^*}}{\sqrt{k-\xi}}, \\
        0;
      \end{array} 
    \right.
\end{equation}
see equation~(\ref{eq:def:H0pm}) in Appendix A for details.

Finally, for $x^*\in W_h$ such that $-h<x_2^*<0$, the associated Green function
$G$ is defined as
\begin{equation}
  \label{eq:G:2}
  G(x;x^*) =\left\{
    \begin{array}{ll}
      G_1(x;x^*), & x \in \mathbb{R}_2^+,\\
      G_2(x;x^*)+G^{\rm in}(x;x^*), & x \in W_h,
    \end{array}
  \right. 
\end{equation}
where
\begin{align*}
  G^{\rm in}(x;x^*) &= \int_{\cal L}\left[  \frac{e^{\bi\mu|x_2-x_2^*|}-e^{\bi\mu(x_2-x_2^*)}}{-4\bi\pi\mu}+\frac{(e^{\bi\mu(2h+x_2^*)}-e^{-\bi\mu x_2^*})(e^{-\bi\mu x_2}-e^{\bi \mu x_2})}{-4\bi\pi(1-e^{2\bi\mu h})\mu}\right]e^{-\bi\xi(x_1-x_1^*)}d\xi
\end{align*}
and $G_1$ and $G_2$ remain the same form as before, but with $H^{\pm}$ redefined
as, for $\xi\in{\cal L}$,
\begin{align}
  \label{eq:H+-:onL2}
  H^{\pm}(\xi;x^*) =& \frac{e^{\bi \xi x_1^*}(e^{\bi \mu(2h+x_2^*)} - e^{-\bi\mu x_2^*})K^-(\xi)}{2(e^{2\bi\mu h}-1)\sqrt{k-\xi}} 
  \pm {\rm p.v.}\int_{\cal L}  \frac{K^-(t)e^{\bi t x_1^*}(e^{\bi \mu(t)(2h+x_2^*)} - e^{-\bi\mu(t) x_2^*})}{2\bi\pi(e^{2\bi\mu(t) h}-1)\sqrt{k-t}(t-\xi)}dt.
\end{align}
It can be seen that $H^{\pm}$ in (\ref{eq:H+-:onL2}) and (\ref{eq:H+-:onL}) have
the same form of limit (\ref{eq:H0pm:L}) as $x_2^*\to 0$.
\subsection{Radiation behavior of $G$ at infinity}
We are concerned with the radiation behavior and far-field pattern of the Green
function $G$. We consider $x_2\geq 0$ first.
\begin{mylemma}
  \label{lem:G1:case1}
  Let $x=(r\cos\alpha,r\sin\alpha)\in \partial D_r$ for $r>0$ and $\alpha = [0,\pi]$. The function $G$ satisfies the following finiteness condition
  \[
    G(x;x^*)= \frac{e^{\bi k
        r}}{\sqrt{r}}(2\pi\bi)^{-1/2}\hf^+(-k\cos\alpha;x^*)k\sin\alpha + {\cal
      O}(r^{-3/2}),\ {\rm as}\ r\to\infty,
  \]
  uniformly for $\alpha\in[0,\pi]$, and
  \begin{equation}
    \label{eq:tan:G1}
    \partial_{\tau(x)} G(x;x^*) = {\cal O}(r^{-3/2}),\ {\rm as}\ r\to\infty,
  \end{equation}
  uniformly for $\alpha\in[0,\pi]$, where ${\tau}(x)=(-\sin\alpha,\cos\alpha)$
  denotes the tangential vector along $\partial D_r$. Moreover, $G(x;x^*)$
  satisfies the Sommerfeld radiation condition
  \begin{equation}
    \label{eq:src:G}
    (\partial_r - \bi k)G(x;x^*)= {\cal O}(r^{-3/2}),\ {\rm as}\ r\to\infty,
  \end{equation}
  uniformly for $\alpha\in[0,\pi]$.
  \begin{proof}
    Without loss of generality, we assume $x_2^*>0$. As $G^{\rm in}(x;x^*)={\cal
      O}(r^{-3/2})$ as $r\to\infty$, we require analyzing $G_1$ only. We
    distinguish three cases:

    (1). $\alpha\in[\frac{3}{4}\pi,\pi]$ so that $x_1<0$, then by Cauchy's
    theorem and by Lemma~\ref{lem:f+}, we get
      \begin{align}
        \label{eq:G1:x1-}
        G_1(x;x^*)=& \frac{1}{2\pi}\int_{ +\infty\bi\to 0\to k} \frac{H^+(\xi;x^*)K^+(\xi)}{2\bi\sqrt{k+\xi}}e^{-\bi \xi x_1+\bi \mu x_2}d\xi \nonumber\\
        &+ \frac{1}{2\pi}\int_{k\to 0\to +\infty\bi} \frac{H^+(\xi;x^*)K^+(\xi)}{2\bi\sqrt{k+\xi}}e^{-\bi \xi x_1-\bi \mu x_2}d\xi.
      \end{align}
      Let $p(\xi) = \frac{1}{2\pi}\hf^+(\xi;x^*)$ and $q^{\pm}(\xi;\alpha) =
      -k\cos\alpha \xi \pm k\sin\alpha \mu(\xi)$. Integration by parts,
      \begin{align*}
        \int_{0}^{+\infty\bi}p(\xi)e^{\bi r q^{\pm}(\xi;\alpha)}d\xi
        =&\left(  \frac{p(\xi)}{\bi r { q^{\pm} }'(\xi;\alpha)}\right)'|_{\xi=0}\frac{e^{\bi r { q^{\pm} }(0;\alpha)}}{\bi r { q^{\pm} }'(0;\alpha)} -\frac{p(0)e^{\bi r { q^{\pm} }(0;\alpha)}}{\bi r {q^{\pm}}'(0;\alpha)}\\
           &+ \int_{0}^{+\infty\bi} \left(  \left( \frac{p(\xi)}{\bi r { q^{\pm}}'(\xi;\alpha)} \right)'\frac{1}{\bi r{ q^{\pm} }'(\xi;\alpha)}\right)'e^{\bi r q^{\pm}(\xi;\alpha)}d\xi,
      \end{align*}
      so that 
      \[
        \left| \int_{0}^{+\infty\bi}p(\xi)e^{\bi r q^{\pm}(\xi;\alpha)}d\xi -
          \frac{p(0)e^{\pm\bi r k\sin\alpha}}{\bi rk \cos\alpha} \right|\leq C r^{-2},
      \]
      where $C$ is independent of $\alpha$. Next, let $\chi$ be a smooth function on
      $\mathbb{R}$ that has a small compact support near $0$ and equals 1 near $0$,
      and we have
      \begin{align*}
        &\int_{0}^{k}p(\xi)e^{\bi r q^+(\xi;\alpha)}d\xi + \int_{k}^{0}p(\xi)e^{\bi r q^-(\xi;\alpha)}d\xi\\
        =&\int_{0}^{\pi}p(k\sin\phi)k\cos\phi  e^{\bi r k\sin(\alpha-\phi)}\chi(\phi-(\alpha-\pi/2)) d\phi \\
        &+ \int_{0}^{\pi}p(k\sin\phi)k\cos\phi e^{\bi r k\sin(\alpha-\phi)} (1-\chi(\phi-(\alpha-\pi/2))) d\phi =:I_1(\alpha) + I_2(\alpha).
      \end{align*}
      By the method of stationary phase \cite[p.334, Prop. 3]{stein93}, we get
      \[
        |I_1(\alpha) - r^{-1/2}(2\pi/\bi)^{1/2}p(-k\cos\alpha)k\sin\alpha e^{\bi rk}| \leq C r^{-3/2}. 
      \]
 On the other hand, integration by parts gives
      \begin{align*}
        \left|I_2(\alpha) - \left[ \frac{p(0)k e^{-\bi rk\sin\alpha}}{-\bi rk\cos\alpha} - \frac{p(0)k e^{\bi rk\sin\alpha}}{-\bi rk\cos\alpha} \right]\right|\leq C r^{-2}.
      \end{align*}
      Combining all the above facts, we get
      \[
        \left| G_1(x;x^*) - r^{-1/2}(2\pi/\bi)^{1/2}p(-k\cos\alpha)k\sin\alpha
          e^{\bi rk}\right| ={\cal O}(r^{-3/2}),
      \]
      uniformly for $\alpha\in[\frac{3\pi}{4},\pi]$. Similarly, one gets
      (\ref{eq:tan:G1}) and the Sommerfeld radiation condition (\ref{eq:src:G})
      uniformly for $\alpha\in[\frac{3\pi}{4},\pi]$.
      
    (2). The second case is $\alpha\in (\frac{\pi}{4},\frac{3\pi}{4})$ so that
      $|\cos\alpha| \leq \sqrt{2}/2\leq \sin\alpha$ and $x_2\geq \sqrt{2}r/{2}$. By
      (\ref{eq:+t-}), Cauchy's theorem and by Lemma~\ref{lem:f+}, we get
      \begin{align*}
        G_1(x;x^*)=&\frac{1}{2\pi}\int_{-\infty}^{+\infty} \hf^+(\xi;x^*)e^{-r\bi\xi \cos\alpha+\bi\mu r\sin\alpha}d\xi.
      \end{align*}
      Setting $p(\xi) = \frac{1}{2\pi}\hf^+(\xi;x^*)$, and $q(\xi) = -\xi\cos\alpha + \mu\sin\alpha$, we get 
      \begin{align*}
        \int_{-\infty}^{-k}p(\xi) e^{\bi r q(\xi)} d\xi
        =&\int_{0}^{+\infty} p(-\sqrt{k^2+t^2})e^{\bi r\sqrt{k^2+t^2}\cos\alpha-tr\sin\alpha}\frac{t}{\sqrt{k^2+t^2}}dt.
      \end{align*}
      As before, integration by parts gives $|\int_{-\infty}^{-k}p(\xi)
        e^{\bi r q(\xi)} d\xi |\leq C r^{-2}$. Similarly, 
      $|\int_{k}^{+\infty}p(\xi) e^{\bi r q(\xi)} d\xi | \leq C r^{-2}$.

      Since
      \[
        \int_{-k}^kp(\xi) e^{\bi r q(\xi)} d\xi = \int_{-\pi/2}^{\pi/2}p(k\sin\phi)e^{\bi r\sin(\alpha-\phi) }\cos\phi d\phi,
      \]
      the method of stationary phase again gives
      \[
        |\int_{-k}^kp(\xi) e^{\bi r q(\xi)} d\xi -
        r^{-1/2}(2\pi/\bi)^{1/2}p(-k\cos\alpha)k\sin\alpha e^{\bi rk}|\leq C r^{-3/2}.
      \]
      Combining all the above facts, we get
      \[
        \left| G_1(x;x^*) - r^{-1/2}(2\pi/\bi)^{1/2}p(-k\cos\alpha)k\sin\alpha e^{\bi rk}\right| \leq C r^{-3/2}.
      \]
      Similarly, one gets (\ref{eq:tan:G1}) and the Sommerfeld radiation
      condition (\ref{eq:src:G}) uniformly for
      $\alpha\in(\frac{\pi}{4},\frac{3\pi}{4})$.

      (3). The last case is $\alpha\in[0,\frac{\pi}{4}]$ so that
      $x_1>0$. By Lemma~\ref{lem:f+} and by Cauchy's theorem, we could rewrite
      $G_1(x;x^*)$ as
      \begin{align}
        \label{eq:G1:x1+}
        G_1(x;x^*) =& \frac{1}{2\pi}\int_{-\infty\bi\to 0 \to -k} \hf_{\rm c}^+(\xi;x^*)e^{-\bi\xi x_1-\bi\mu x_2}d\xi + \frac{1}{2\pi}\int_{-k\to 0\to -\infty\bi}\hf^+(\xi;x^*)e^{-\bi\xi x_1+\bi\mu x_2}d\xi,
      \end{align}
      where since branch cut has been crossed, we have defined
      \[
        \hf_c^+(\xi;x^*) = \frac{e^{-\bi \mu x_2^*}(1-e^{-2\bi\mu h})}{-2\bi\mu}e^{\bi\xi x_1^*}+\frac{H^-(\xi;x^*)(1-e^{-2\bi\mu h})}{2\bi\sqrt{k+\xi}K^-(\xi)}.
      \]
      Then, as in case (1), we get 
      \[
        \left| G_1(x;x^*) -  r^{-1/2}(2\pi\bi)^{-1/2}\hf^+(-k\cos\alpha;x^*)k\sin\alpha
          e^{\bi k r}\right|\leq C r^{-3/2},
      \]
      (\ref{eq:tan:G1}) and the Sommerfeld radiation condition (\ref{eq:src:G})
      uniformly for $\alpha\in[0,\frac{\pi}{4}]$.

    The lemma follows from combining the above three cases. 
  \end{proof}
\end{mylemma}
As for $x_2\in[-h,0]$, we need to distinguish two situations: $x_1>0$ and
$x_1<0$. The following lemma reveals the asymptotic behavior of $G$ for $x$
outside the semi-waveguide, i.e., $x_1>0$ and $x_2\in[-h,0]$.
\begin{mylemma}
  \label{lem:G2:case1}
For $x_1>0$, the Green function $G(x;x^*)$ in $\overline{W_h}$ satisfies the following Sommerfeld radiation condition 
  \[
    (\partial_{x_1} - \bi k) G(x;x^*) = {\cal O}(x_1^{-3/2}),\ {\rm as}\ x_1\to+\infty,
  \]
  uniformly for $x_2\in[-h,0]$. Moreover, $G$ satisfies the following
  finiteness condition 
  \[
    G(x;x^*)= \frac{e^{\bi k x_1}}{\sqrt{x_1}}{\cal O}(x_1^{-1}),\ {\rm as}\ x_1\to+\infty,
  \]
  uniformly for $x_2\in[-h,0]$.
  \begin{proof}
    We consider case $x_2^*>0$ only.
    For $x_1>0$, we could use Cauchy's theorem to rewrite $G=G_2$ in
    (\ref{eq:sol:G2}) by (\ref{eq:+t-}) as follows,
    \begin{align*}
      G(x;x^*) =\Phi^k(x';x^*) - \Phi^k(x;x^*) 
      -\int_{-k\to 0\to-\infty\bi}\frac{e^{-\bi\xi x_1}H^-(\xi)(\cos(\mu x_2)-\cos(\mu(x_2+2h)))}{2\pi\bi K^-(\xi)\sqrt{k+\xi}}d\xi,
    \end{align*}
    where $x'=(x_1,-x_2-2h)$. The results then follow from similar arguments in
    the proof of Lemma~\ref{lem:G1:case1} and from the method of stationary
    phase as in \cite[p.334, Prop. 3]{stein93}.
  \end{proof}
\end{mylemma}
When $x_1<0$, $G$ contains propagating eigensolutions $\{\phi_m = e^{-\bi \xi_m
  x_1}\sin(\mu_m x_2)\}_{m=1}^M$, with $M$ chosen as the largest integer such
that $\mu_m=\frac{m\pi}{h}\in [0,k]$ and $\xi_m=\sqrt{k^2-\mu_m^2}\geq 0$ for
all $m=1,\cdots M$, as was done in \cite{tan99}. In fact, $G$ satisfies the MRC
condition \cite{bongouhaz11,bontil01a, bontil01b} in $W_h$ as $x_1\to-\infty$,
as shown below.
\begin{mylemma}
  \label{lem:G2:case2}
For $x_1<0$, the Green function $G(x;x^*)$ in $\overline{W_h}$ can be expanded as 
\begin{equation}
  \label{eq:exp:G}
  G(x;x^*) = \sum_{m=1}^{M}c_m(x^*)\phi_m(x) + G^{\rm dec}(x;x^*),
\end{equation}
where the first term represents the guided part with known Fourier coefficients
$c_m(x^*)$, and the second term $G^{\rm dec}$ represents the super-algebraically
decaying part satisfying $G^{\rm dec}(x;x^*)= {\cal O}(|x_1|^{-N})$,\ {\rm as}\
$x_1\to-\infty$, uniformly for $x_2\in[-h,0]$, for any integer $N>0$.
  \begin{proof}
    We prove case $x_2^*>0$ only. For $x_1<0$, we could use residual theorem to rewrite $G=G_2$ in
    (\ref{eq:sol:G2}) by (\ref{eq:+t-}) as follows: 
    \begin{align*}
      G(x;x^*) =&\sum_{m=1}^{M}c_m(x^*)\phi_m(x)+\frac{1}{2\pi}\int_{-\infty+\delta_0\bi}^{+\infty+\delta_0\bi}\hf^+(\xi;x^*)\frac{e^{-\bi\mu x_2}-e^{\bi\mu (x_2+2h) }}{1-e^{2\bi\mu h}} e^{-\bi \xi x_1}d\xi 
    \end{align*}
    where $c_m(x^*)=\frac{\mu_m}{\bi h\xi_m}\hf^+(\xi_m;x^*)$ for $1\leq m\leq
    M-1$, 
    \[
      c_M(x^*)=\left\{
        \begin{array}{lc}
          \frac{\mu_M}{\bi h\xi_M}\hf^+(\xi_M;x^*) & \xi_M>0,\\
          \left.\frac{d}{d\xi}\left( \frac{\hf^+(\xi;x^*)\sin(\mu
          x_2)e^{-\bi\mu h}\xi^2}{\bi\sin(\mu h)} \right)\right|_{\xi=0} & \xi_M=0,
        \end{array}
      \right.
\]
and the positive constant $\delta_0$ is chosen such that $1-e^{2\bi\mu h}\neq 0$
for all $\xi\in(0,\delta_0)$; we note that the Green function in \cite{tan99}
is invalid when $\xi_M=0$ since it corresponds to a pole of order $2$, but
not $1$. The results then follow immediately.
  \end{proof}
\end{mylemma}

%

\section{Sommerfeld radiation condition and wellposedness}

\subsection{Radiating solution and Far-field pattern}
Throughout this section, let $\Omega_R^+=\Omega\backslash\overline{D_R}$ and
$\Omega_R^-=\Omega\cap D_R$ be the exterior and interior region separated by
$\partial D_R$, respectively where we recall $D_R$ is the disk of radius $R>0$
centered at the origin. Let $\Gamma_R=\partial D_R\cap \Omega$ be the open arc
of $\partial D_R$, $\Gamma_R^{\rm
  ext}=\partial\Omega_R^+\backslash\overline{\Gamma_R}$ and $\Gamma_R^{\rm
  int}=\partial\Omega_R^-\backslash\overline{\Gamma_R}$ be the exterior and
interior part of $\partial\Omega$ separated by the arc $\Gamma_R$.
We first give the definition of a radiating solution as
follows.
\begin{mydef}
  A solution $u$ of the problem (\ref{eq:gov:tot1}) and (\ref{eq:gov:tot2}) is
  called radiating if for some $R>0$, $u\in C^2(\overline{\Omega_R^+})$
  satisfies the following universal-direction Sommerfeld radiation
  condition (uSRC): for $x=(r\cos\alpha,r\sin\alpha)\in
  \overline{\Omega_R^+\cap\mathbb{R}^2_+}$,
  \begin{equation}
    \label{eq:src:1}
    (\partial_r- \bi k) u(x) = {\cal O}(r^{-3/2}),\ {\rm as}\ r\to\infty,
  \end{equation}
  uniformly for all $\alpha\in[0,\pi]$, and for
  $x=(x_1,x_2)\in\overline{\Omega_R^+\cap W_h}$,
  \begin{equation}
    \label{eq:src:2}
  (\partial_{x_1} - \bi k)u(x) = {\cal O}(x_1^{-3/2}),\ {\rm as}\ x_1\to+\infty,
  \end{equation}
  uniformly for all $x_2\in[-h,0]$.
\end{mydef}
In comparison to the SRC conditions proposed in \cite{luhu19,hulurat20}, we
don't require $\sup_{x\in \Omega_R^+}|x|^{1/2}|u(x)|<+\infty$, but we shall see
below that this condition is automatically satisfied. Clearly, $G$ is radiating
and satisfies the uSRC. Moreover, any radiating solution satisfies the
following Green's representation formula.
\begin{mytheorem}
  \label{thm:green}
  Then any radiating
  solution $u$ vanishing on $\partial\Omega$  satisfies
  \begin{equation}
    \label{eq:greenformula}
    u(x) = \int_{\Gamma_R}\left[ u(x')\partial_{\nu'}G(x;x') - \partial_{\nu'} u(x') G(x;x') \right]ds(x'),\quad {\rm for}\quad x\in\Omega_R^+. 
  \end{equation}
  \begin{proof}
    The proof is analogous to the proof of Theorem 2.5 in \cite{colkre13}.
    It can be shown from the proof that
    \[
      \int_{\partial B(0,r)\cap \Omega_R^+} |u(x)|^2ds(x) = {\cal O}(1),\quad
      r\to\infty,
    \]
    which has been proposed explicitly as part of a radiation condition in
    \cite{luhu19,hulurat20}. 
  \end{proof}
\end{mytheorem}
Based on the Green's formula (\ref{eq:greenformula}) and the far-field
pattern of $G$, we get the finiteness condition and far-field pattern of a
radiating wave $u$.
\begin{mylemma}
  \label{lem:farfield}
  Suppose $u$ denotes a radiating solution vanishing on
  $\partial\Omega$. For $x = (r\cos\alpha,r\sin\alpha)\in \overline{\Omega\cap \mathbb{R}_2^+}$ with
  $\alpha\in[0,\pi]$, $u$ has the asymptotic behavior of an outgoing wave
  \begin{equation}
    \label{eq:u}
    u(x) = \frac{e^{\bi k r}}{\sqrt{r}}\left[ u_{\infty}(\hat{x}) + {\cal O}(\frac{1}{r}) \right],\quad r\to\infty,
  \end{equation}
  uniformly for $\hat{x}=(\cos\alpha,\sin\alpha)$ and any $\alpha\in[0,\pi]$.
  Here, $u_{\infty}(\hat{x})$ denotes the half-plane far-field pattern of $u$
  satisfying
  \begin{equation}
    \label{eq:farfield}
    u_\infty(\hat{x}) = (2\pi\bi)^{-1/2}k\sin\alpha \int_{\Gamma_R}\left[u(x')\partial_{\nu'}\hat{f}^+(-k\cos\alpha;x') - \partial_{\nu'} u(x') \hat{f}^+(-k\cos\alpha;x') \right] ds(x').
  \end{equation}
  If $x = (x_1,x_2)\in \overline{\Omega\cap W_h}$, then any radiating solution $u$ satisfies 
  \[
    u(x) = \frac{e^{\bi k x_1}}{\sqrt{x_1}}{\cal O}(\frac{1}{x_1}),\quad x_1\to+\infty,
  \]
  uniformly for all $x_2\in[-h,0]$; i.e., the far-pattern along $x_1$-direction
  is $0$.
\end{mylemma}

\subsection{Implicit transparent boundary condition}
We now derive an implicit transparent boundary condition (ITBC) on $\Gamma_R$ for
any radiating solution $u\in C^{2}(\overline{\Omega_R^+})$ vanishing on
$\partial\Omega$. Before proceeding, we first introduce the following four
integral operators on the open arc $\Gamma_R$,
\begin{align*}
  [{\cal S} \phi](x) &= 2\int_{\Gamma_R} G(x;y) \phi(y) ds(y),\quad
  [{\cal K} \phi](x) = 2{\rm p.v.}\int_{\Gamma_R} \partial_{\nu(y)}G(x;y) \phi(y) ds(y),\\
  [{\cal K}' \phi](x) &= 2{\rm p.v.}\int_{\Gamma_R} \partial_{\nu(x)}G(x;y) \phi(y) ds(y),\quad
  [{\cal T} \phi](x) = 2{\rm f.p.}\int_{\Gamma_R} \partial^2_{\nu(x)\nu(y)}G(x;y) \phi(y) ds(y),
\end{align*}
where f.p. indicates the finite part integral. Since $G(x;y)-\Phi(x;y)\in
C^{\infty}(\Omega_R^+)$, they satisfy according to \cite[Thm. 7.1]{mcl00} the
following classic mapping properties $\cS: C^{\infty}_{\rm comp}(\Gamma_R)\to
C_0^{\infty}(\overline{\Gamma_R})$, $\cK: C_{\rm comp}^{\infty}(\Gamma_R)\to
C_0^{\infty}(\overline{\Gamma_R})$, $\cKd: C_{\rm comp}^{\infty}(\Gamma_R)\to
C^{\infty}(\overline{\Gamma_R})$, and $\cT: C_{\rm comp}^{\infty}(\Gamma_R)\to
C^{\infty}(\overline{\Gamma_R})$. Here, the subscript $0$ indicates that the
function vanishes at the endpoints of $\Gamma_R$ since
$\partial_{\nu(y)}G(x;y)=G(x;y)=0$, for any $x\in\Gamma_R^{\rm ext}$ and any
$y\in\Gamma_R$. First, we extend the mapping definitions of the four integral
operators in standard Sobolev spaces.
\begin{mylemma}
  \label{lem:ext:SK} 
  We can uniquely extend the operator $\cS$ as a bounded operator from
  $H^{-1/2}(\Gamma_R)\to\wtd{H^{1/2}}(\Gamma_R)$, the operator $\cK$ as a
  compact (and certainly bounded) operator from
  $\wtd{H^{1/2}}(\Gamma_R)\to\wtd{H^{1/2}}(\Gamma_R)$, the operator $\cKd$ as a
  compact operator $H^{-1/2}(\Gamma_R)\to H^{-1/2}(\Gamma_R)$, and the operator
  $\cT$ as a bounded operator $\wtd{H^{1/2}}(\Gamma_R)\to H^{-1/2}(\Gamma_R)$.
  Moreover, we have the decomposition $\cS = \cS_p + {\cal L}_p$ such that
  $\cS_p: H^{-1/2}(\Gamma_R)\to \wtd{H^{1/2}}(\Gamma_R)$ is positive and bounded
  below, i.e., for some constant $c>0$,
  \[
    {\rm Re}( <{\cal S}_p\phi,\phi>_{\Gamma_R} )\geq c ||\phi||^2_{H^{-1/2}(\Gamma_R)},
  \]
  for any $\phi\in H^{-1/2}(\Gamma_R)$, and ${\cal L}_p:H^{-1/2}(\Gamma_R)\to
  \wtd{H^{1/2}}(\Gamma_R)$ is compact.
  \begin{proof}
    Let ${\cal S}_R: H^{-1/2}(\partial\Omega_R^-)\to H^{1/2}(\partial\Omega_R^-)$ be
    the free-space single-layer operator \cite[Thm. 7.1]{mcl00}
    \[
      [{\cal S}_R\phi](x) = 2\int_{\partial\Omega_R^-} \Phi^k(x;y)\phi(y)ds(y).
    \]
    According to \cite[Thm. 7.6]{mcl00},
    we can decompose ${\cal S}_R = {\cal S}_{R,p} + {\cal L}_R$ such
    that for any $\phi\in H^{-1/2}(\partial\Omega_R^-)$,
    \[
      {\rm Re}( <{\cal S}_{R,p}\phi,\phi>_{\partial\Omega_R^-} )\geq c||\phi||_{H^{-1/2}(\partial\Omega_R^-)},
    \]
    and ${\cal L}_R$ is compact. 

    Let $G_2(x;y) = G(x;y)-\Phi(x;y)\in C^{\infty}(\overline{\Omega\cap\Omega_h})$ for
    any $y\in \Gamma_R$. We notice that $G_2$ may not be well-defined in domain
    $\Omega$ since $\Omega$ may not be a subset of $\Omega_h$. However, we could
    extend $G_2$ to $\mathbb{R}^2$ to make itself at least $C^2(\mathbb{R}^2)$
    \cite[Thm. A.4]{mcl00}. The associated single-layer potential with kernel
    $G_2$ then defines a compact operator $\cS_2:
    H^{-1/2}(\partial\Omega_R^-)\to H^{1/2}(\partial\Omega_R^-)$. Now, we define for any $x\in \partial\Omega_R^-$
    and any $\phi\in C_{\rm comp}^{\infty}(\Gamma_R)$ that
    \begin{align}
      \label{eq:ext:cS}
      [\cS\phi](x)&:=[\chi_R\cS_2(\chi_RE{\phi})](x) + [\chi_R{\cal S}_R(\chi_RE{\phi})](x)\nonumber\\
                  &=\left\{  [\chi_R\cS_2(\chi_RE{\phi})](x) + [\chi_R{\cal L}_R(\chi_RE{\phi})](x)\right\} + [\chi_R\cS_{R,p}(\chi_RE{\phi})](x), 
    \end{align}
    where we recall that $E: H^{-1/2}(\Gamma_R)\to H^{-1/2}(\partial\Omega_R^-)$
    is the bounded extension operator such that $(E{\phi})|_{\Gamma_R}=\phi$,
    $\chi_R\in C_{\rm comp}^{\infty}(\mathbb{R}^2)$ is chosen such that
    $\chi_R(x)=1$ for all $x\in \Gamma_R$ and has a support in a sufficiently
    small neighborhood of $\Gamma_R$. One
    easily verifies that $[\cS\phi](x)$ remains as before for $x\in \Gamma_R$,
    but becomes zero elsewhere on $\Gamma_R^{\rm int}$. Thus,
    \[
      ||\cS \phi||_{\wtd{H^{-1/2}}(\Gamma_R)}\leq C
      (||\cS_2|| + ||\cS_R||)||E||\cdot ||\phi||_{H^{-1/2}(\Gamma_R)}.
    \]
    Since $C^{\infty}_{\rm comp}(\Gamma_R)$ is dense in $H^{-1/2}(\Gamma_R)$,
    ${\cal S}$ defined in (\ref{eq:ext:cS}) can be uniquely extended as a bounded
    operator from $H^{-1/2}(\Gamma_R)$ to $\wtd{H^{1/2}}(\Gamma_R)$. Now, define
    for any $\phi\in H^{-1/2}(\Gamma_R^-)$,
    \[
      {\cal L}_p{\phi}= E^*\chi_R\cS_2(\chi_RE{\phi}) + E^*\chi_R{\cal L}_R(\chi_RE{\phi}),\quad {\cal S}_p{\phi} = E^*\chi_R{\cal S}_{R,p}(\chi_RE\phi),
    \]
    where $E^*: H^{1/2}(\partial\Omega_R^-)\to \wtd{H^{1/2}}(\Gamma_R)$ is the
    ajoint operator of $E$. Since for any $\phi,\psi\in C_{\rm comp}^{\infty}(\Gamma_R)$,
    \[
      <E^*\cS\phi, \psi>_{\Gamma_R} = <\cS\phi, E\psi>_{\partial\Omega_R^-} =
      <\cS\phi, E\psi|_{\Gamma_R}>_{\Gamma_R} = (\cS\phi,\psi)_{\Gamma_R}=<\cS\phi,\psi>_{\Gamma_R},
    \]
    we have $\cS = E^*\cS = {\cal L}_p + {\cal S}_p$.
    Thus,
    \[
      {\rm Re}(<{\cal S}_p\phi,\phi>_{\Gamma_R})={\rm Re}(<{\cal S}_{R,p}\chi_RE\phi,\chi_RE\phi>_{\partial\Omega_R^-}) \geq
      c||\chi_RE{\phi}||_{H^{-1/2}(\partial\Omega_R^-)}\geq
      c||\phi||_{H^{-1/2}(\Gamma_R)},
    \]
    where the last inequality is due to $(\chi_RE\phi)|_{\Gamma_R}=\phi$. The
    compactness of ${\cal L}_p$ follows immediately from the compactness of $\cS_2$
    and ${\cal L}_R$.
    Similarly, we define for any $\phi\in \wtd{H^{1/2}}(\Gamma_R)$ that 
    \begin{align*}
      {\cal K}\phi= \chi_R{\cal K}_2\phi + \chi_R{\cal K}_R\phi = E^*\chi_R{\cal K}_2\phi + E^*\chi_R{\cal K}_R\phi,
    \end{align*}
    where ${\cal K}_2$ and ${\cal K}_R$ are double-layer operators of smooth kernel
    $\partial_{\nu(y)}G_2(x;y)$ and the free-space double-layer kernel
    $\partial_{\nu(y)}\Phi(x;y)$ over the smooth boundary $\partial B(0,R)$,
    respectively. It is clear that ${\cal K}\phi\in C_0^{\infty}(\Gamma_R)$ remains the same
    when $\phi\in C_{\rm comp}^{\infty}(\Gamma_R)$. The compactness of ${\cal K}$ then
    follows immediately from the compactness of ${\cal K}_2$ and ${\cal K}_R$. The
    operator $\cKd$ is compact since it is the dual operator of ${\cal K}$. The
    boundedness of ${\cal T}$ is similar to prove.
  \end{proof}
\end{mylemma}

In previous work \cite{hulurat20, luhu19, baohuyin18}, the direct boundary
integral equation $u - {\cal K}u =-\cS\partial_{\nu} u$ for wave field $u$ and
its normal derivative $\partial_{\nu} u$ was adopted as the TBC condition. A
potential difficulty is that one should carefully choose the open arc $\Gamma_R$
to ensure $-k^2$ is not a resonant frequency of the Laplace equation for the
interior domain $\Omega_R^-$. To resolve this issue, we shall propose an
indirect version of TBC condition in the following. As inspired in
\cite{colkre13, kre14}, we introduce an auxiliary density function $\phi\in
H^{-1/2}(\Gamma_R)$ to represent $u$ in $\Omega_R^+$ as
\begin{equation}
  \label{eq:rep}
  u(x) = \int_{\Gamma_R}\left[ G(x;y)\phi(y) + \bi\eta
    \partial_{\nu(y)}G(x;y)(E^*[\cS_0^2E\phi](y)) \right] ds(y)\in H^{1}_{\rm loc}(\Omega_R^+),
\end{equation}
for some positive constant $\eta$, where $\cS_0:H^{-1/2}(\partial D_R)\to H^{1/2}(\partial
D_R)$ is the modified single-layer potential operator of two-dimensional
Laplace equation as introduced in \cite[p. 134\&169]{kre14}. By the standard
jump conditions and by Lemma~\ref{lem:ext:SK}, we can verify that on
$\Gamma_R$,
\begin{align}
  \label{eq:tbc:1}
  2\gamma u &= \cS \phi + \bi \eta \cK(E^*\cS_0^2 E\phi) + \bi \eta E^*\cS_0^2E\phi \in \wtd{H^{1/2}}(\Gamma_R),\\
  \label{eq:tbc:2}
  2\partial_{\nu} u &= \cKd \phi - \phi + \bi \eta \cT(E^*\cS_0^2 E\phi) \in H^{-1/2}(\Gamma_R),
\end{align}
for any $\phi\in H^{-1/2}(\Gamma_R)$. The two equations (\ref{eq:tbc:1}) and
(\ref{eq:tbc:2}) linking $u$ and $\partial_{\nu}u$ on $\Gamma_R$ together via an
unknown density function $\phi$ is called the implicit TBC (ITBC) condition in
the following.

We are in the position to propose sharper radiation conditions and prove the
well-posedness for the scattering problem by each of the following two incident
waves: (1) cylindrical incident wave; (2) plane incident wave. 

\subsection{cylindrical incident wave}
Suppose $u^{\rm inc}(x)= \frac{\bi}{4}H_0^{(1)}(k|x-x^*|)$ with $x^*\in\Omega$.
As inspired by \cite{luhu19,hulurat20}, we directly enforce $u^{\rm tot}$ the
uSRC condition (\ref{eq:src:1}) and (\ref{eq:src:2}), and obtain the following
boundary value problem: Find the radiating wave $u^{\rm tot}=u^{\rm inc} +
u^{\rm sc}$ with $u^{\rm sc}\in H^{1}_{\rm loc}(\Omega)\cap
C^2(\overline{\Omega_R^+})$, such that $u^{\rm tot}$ solves (\ref{eq:gov:tot1}),
with $-\delta(x-x^*)$ in place of the r.h.s, and (\ref{eq:gov:tot2}). In fact,
$u^{\rm tot}$ represents the Green function for scattering surface
$\partial\Omega$ excited by source $x^*$.

Let $\tilde{u}^{\rm inc}=(1-\chi_0)u^{\rm inc}\in H_{\rm loc}^{1}(\Omega)$ where
$\chi_0\in C_{\rm comp}^{\infty}(\mathbb{R}^2)$ satisfies $\chi_0=1$ in a
sufficiently small neighborhood of the source $x^*$ and has a sufficiently small
support; note that $u^{\rm inc}\notin H_{\rm loc}^{1}(\Omega)$. Thus, for
$u=u^{\rm sc}+ \tilde{u}^{\rm inc}\in H_{\rm loc}^{1}(\Omega)$ which is exactly
$u^{\rm tot}$ in $\Omega_R^+$ for some $R>0$ and hence radiating at infinity.
Let $H_{\rm int}^1(\Omega_R^-)=\{u\in H^1(\Omega_R^-): u=0\ {\rm on}\
\Gamma_R^{\rm int}\}$, and $V=H_{\rm int}^{1}(\Omega_R^-)\times
H^{-1/2}(\Gamma_R)$ equipped with the natural product norm. Consequently, we
seek $(u,\phi)\in V$ such that, in a distributional sense,
\begin{align}
  \label{eq:u1:1}
  \Delta u +  k^2 u &= f,\quad{\rm on}\quad\Omega_R^-,\\
  \label{eq:u2:1}
  \gamma u &= 0,\quad{\rm on}\quad \Gamma_R^{\rm int},\\
  \label{eq:u3:1}
  (u,\phi)\ {\rm satisfies}\ &{\rm ITBC}\ (\ref{eq:tbc:1})\ {\rm and}\ (\ref{eq:tbc:2})\ {\rm on}\ \Gamma_R.
\end{align}
Here we note that $f=\Delta \tilde{u}^{\rm inc}+k^2\tilde{u}^{\rm inc}\in
\wtd{H^{-1}}(\Omega_R^-)$ since $\Delta \tilde{u}^{\rm inc} + k^2 \tilde{u}^{\rm
  inc} = 0$ on $\Omega\backslash\overline{\Omega_R^-}$, and $\gamma u \in
\wtd{H^{1/2}}(\Gamma_R)$ due to (\ref{eq:u2:1}).

Thus, an equivalent variational formulation can be posed as: Let $V=
H_{\rm int}^{1}(\Omega_R^-)\times H^{-1/2}(\Gamma_R)$, and $a: V\times V\to \mathbb{C}$
be the following bounded sesquilinear form,
\begin{align}
  a((u,\phi),(v,\psi)) :=& (\nabla u,\nabla v)_{\Omega_R^-} - k^2(u,v)_{\Omega_R^-} - \frac{1}{2}<\cKd \phi - \phi + \bi \eta \cT(E^*\cS_0^2 E\phi), \gamma v>_{\Gamma_R} \nonumber\\
                         &+ <[\cS \phi + \bi \eta \cK(E^*\cS_0^2 E\phi) + \bi \eta E^*\cS_0^2E\phi]/4-\gamma u/2, \psi>_{\Gamma_R}.
\end{align}
Find $(u,\phi)\in V$ such that for any $(v,\psi)\in V$,
\begin{equation}
  \noindent{\rm (P1):}\quad\quad a((u,\phi),(v,\psi))=<f,v>_{\Omega_R^-}.
\end{equation}
We have the following well-posedness result.
\begin{mytheorem}
  For any incident cylindrical wave $u^{\rm inc}(x)=\frac{\bi}{4}H_0^{(1)}(k|x -
  x^*|)$ with $x^*\in \Omega$ and any $k>0$, there exists a unique radiating
  solution $u^{\rm tot}= u^{\rm sc}+u^{\rm inc}$ with $u^{\rm sc}\in H^1_{\rm
    loc}(\Omega)$.
  \begin{proof}
    According to the definition and Lemma~\ref{lem:ext:SK}, we can decompose
    $a=a_1+a_2$ where
    \begin{align}
      a_1((u,\phi),(v,\psi)) =& (\nabla u,\nabla v)_{\Omega_R^-} - k^2(u,v)_{\Omega_R^-} + <\phi, \gamma v>_{\Gamma_R}/2 - <\gamma u, \psi>_{\Gamma_R}/2  \nonumber\\
                              &+ <{\cal S}_p \phi, \psi>_{\Gamma_R}/4 + <\bi \eta E^*\cS_0^2E\phi,\psi>_{\Gamma_R}/4,\\
      a_2((u,\phi),(v,\psi)) =& - <{\cal K}'\phi + \bi \eta \cT(E^*\cS_0^2E\phi),\gamma v>_{\Gamma_R}/2 + <{\cal L}_p \phi + \bi \eta \cK(E^*\cS_0^2 E\phi), \psi>_{\Gamma_R}/4.
    \end{align}
    According to Lemma~\ref{lem:ext:SK}, $a_1$ is coercive on $V$ as
    \[
      {\rm Re}(a_1((u,\phi),(u,\phi)))\geq ||u||_{H^1(\Omega_R^-)}^2 -
      C||u||_{L^2(\Omega_R^-)}^2 + c||\phi||_{H^{-1/2}(\Gamma_R)}^2,
    \]
    and the bounded linear operator associated with $a_2$ is compact.
    Consequently, $a$ is Fredholm of index zero \cite[Thm. 2.34]{mcl00}. Next,
    we prove the uniqueness. Suppose there exists $(u,\phi)\in V$ such that
    \[
      a((u,\phi),(v,\psi)) = 0,\quad \forall(v,\psi)\in V.
    \]
    Then, $(u,\phi)\in V$ solves the problem (\ref{eq:u1:1}-\ref{eq:u3:1}) in
    $\Omega_R^-$ with $0$ in place of $f$. However, we can directly extend $u$
    to $\Omega_R^+$ by (\ref{eq:rep}), and denote the wave field by $u^+$.
    The jump conditions and the ITBC condition (\ref{eq:tbc:1}) and
    (\ref{eq:tbc:2}) imply that $\gamma u^+=\gamma u$ and $\partial_{\nu}
    u^+=\partial_{\nu} u$ on $\Gamma_R$, so that
    \[
      \tilde{u} = \left\{
        \begin{array}{lc}
          u^+,&{\rm on}\quad \Omega_R^+,\\
          u,&{\rm on}\quad \Omega_R^-,
        \end{array}
      \right.
    \]
    in $H^1_{\rm loc}(\Omega)$ solves the homogeneous scattering problem
    (\ref{eq:gov:tot1}) and (\ref{eq:gov:tot2}) and is radiating at infinity.
    Since the uSRC condition automatically satisfies the UPRC/ASR condition, we
    get $\tilde{u}=0$ on $\Omega$ according to the uniqueness result in
    \cite[Thm. 4.1]{chaels10}.

    Now we show $\phi=0$ as well. To this purpose, we use (\ref{eq:rep})
    to define a solution $u_-\in H_{\rm loc}^1(\Omega_b^-)$, where
    $\Omega_b^-=\Omega_b\backslash\overline{\Omega_R^+}$ and we recall that
    $\Omega_b$ is the background domain used before. Then, the jump condition
    and $u_+=0$ in $\Omega_R^+$ gives
    \[
      u_-= - \bi\eta E^*{\cal S}_0^2E\phi\in \wtd{H^{1/2}}(\Gamma_R),\quad
      \partial_{\nu} u_-= \phi\in H^{-1/2}(\Gamma_R). 
    \]
    On the other hand, from the representation (\ref{eq:rep}) and the expansion
    (\ref{eq:exp:G}) of the Green function $G$, $u_-(x)$ admits the
    following unique expansion
    \[
      u_-(x) = \sum_{m=1}^M \hat{u}_{-,m}\phi_m(x) + u^{\rm dec}(x),
    \]
    in the waveguide $W_h$ for $x_1<0$ and $x_2\in(-h,0)$ with Fourier
    coefficients $\hat{u}_{-,m}$ and $u^{\rm dec}(x)={\cal O}(|x_1|^{-N})$ for
    any $N\geq 0$, such that for any $H>0$, in the domain $W_H$ bounded by
    $\Gamma_h$, $\Gamma_R$, the horizontal axis, and the vertical line
    $V_H=\{(-H,y): 0<y<h\}$, we have by Green's identity and the expansion
    of $u_-$, and the orthogonality of $\{\phi_m\}_{m=1}^{M}$ along $V_H$, that
    \begin{align*}
      (\nabla u_-,\nabla u_-)_{W_H} - k^2(u_-,u_-)_{W_H} =& <\partial_{\nu} u_-,u_->_{\Gamma_R\cup V_H} \\
      =& \bi\eta||{\cal S}_0E\phi||_{L^2(\partial D_R)} + \frac{h\bi}{2}\sum_{m=1}^{M}\xi_m|\hat{u}_m|^2 + {\cal O}(H^{-N}). 
    \end{align*}
    Considering the imaginary part of the above and letting $H\to \infty$, we
    get ${\cal S}_0E\phi=0$ in the $L^2(\partial D_R)$ so that $E\phi=0$ in
    $H^{-1/2}(\partial D_R)$ due to the bijectivity of $\cS_0$ \cite[p.
    169]{kre14}. Consequently, $\phi=E\phi|_{\Gamma_R}=0$ in
    $H^{-1/2}(\Gamma_R)$, and the proof is concluded since $(f,v)_{\Omega_R^-}$
    defines a bounded linear operator in $V^*$.
  \end{proof}
\end{mytheorem}

\subsection{Plane incident wave}
Now let $u^{\rm inc}(x) = e^{\bi k(\cos\theta x_1 - \sin\theta x_2)}$, where
$\theta\in(0,\pi)$ denotes the angle between the incident direction and the
positive horizontal axis. Unlike the previous case, we could no longer enforce
$u^{\rm tot}$ the uSRC condition as it contains $u^{\rm inc}$ and two reflected
plane waves of different phases,
\[
  u^{\rm ref}_-(x) = -e^{\bi k(\cos\theta x_1 +\sin\theta x_2)},\ {\rm and}\ u^{\rm ref}_+(x) = -e^{2k\sin\theta h}e^{\bi k(\cos\theta x_1 +\sin\theta x_2)},
\]
due to the two horizontal parts of $\partial\Omega$ of different heights at
negative and positive infinity, respectively.

Let $L_\theta=\{(t\cos\theta,t\sin\theta): t>0\}$ and let
$\Omega_{b,\theta}^{\mp}= \{x=(x_1,x_2)\in\Omega_{b}:
\mp(x_2\cos\theta-x_1\sin\theta)>0\}$ be the two sub-domains of $\Omega_b$ on
the left and right of $L_\theta$, respectively. In the following, we show that
$u^{\rm tot}$, after subtracting a background solution in $\Omega_b$, satisfies
the uSRC condition (\ref{eq:src:1}) and (\ref{eq:src:2}), i.e., it is radiating.
\begin{mylemma}
  \label{lem:backsol}
  The following function
  \begin{equation}
    \label{eq:def:backplane}
    u^{\rm tot}_b(x) = u^{\rm inc}(x)+ \int_{L_\theta}\partial_{\nu(y)}G(x;y)(u^{\rm ref}_+(y) - u^{\rm ref}_-(y))ds(y)  +        u^{\rm ref}_\mp(x), \quad x\in \Omega_{b,\theta}^\mp,
  \end{equation}
  can be continuously extended to be $C^\infty(\Omega_b)$, and vanishes on $\partial\Omega_h$, where ray
   and $\nu(y)$ is the normal
  vector along ray $L_\theta$ towards domain $\Omega_{b,\theta}^+$.
  \begin{proof}
    Notice that the integral in (\ref{eq:def:backplane}) is related to a
    double-layer potential function and its integrand decays of rate
    $|y|^{-3/2}$ as $|y|\to\infty$ due to Lemma~\ref{lem:G1:case1}, as $\nu(y)$
    in fact is the tangential vector of a unit circle. Thus, for any $x^*\in
    L_\theta$, by the standard jumping conditions,
    $u^{\rm tot}$ can be extended to be continuous across $L_\theta$ since
    \begin{align*}
      \lim_{x\to x^{*,\pm}}u^{\rm tot}_b(x)
                                       &=u^{\rm inc}+{\rm p.v.}\int_{L_\theta}\partial_{\nu(y)}G(x^*;y)(u^{\rm ref}_+(y) - u^{\rm ref}_-(y))ds(y) +\frac{1}{2}(u^{\rm ref}_+(y) + u^{\rm ref}_-(y)),
    \end{align*}
    and its normal derivative is continuous across $L_\theta$ since
    \begin{align*}
      \lim_{h\to 0^\pm}\partial_{\nu(x^*)}u^{\rm tot}_b(x^*+h\nu(x^*))
                                                                    &=\partial_{\nu(x^*)}u^{\rm inc}(x^*)+{\rm f.p.}\partial_{\nu(x^*)}\int_{L_\theta}\partial_{\nu(y)}G(x^*;y)(u^{\rm ref}_+(y) - u^{\rm ref}_-(y))ds(y), 
    \end{align*}
    where we notice that $\partial_{\nu(x^*)}u^{\rm ref}(x^*)=0$. The
    $C^\infty-$smoothness at any point $x^*\in L_\theta$ then follows
    immediately from applying Green's representation formula in a neighborhood
    of $x^*$.
  \end{proof}
\end{mylemma}
Thus, we can pose the following boundary value problem: Find $u^{\rm tot}\in
H^{1}_{\rm loc}(\Omega)\cap C^2(\overline{\Omega_R^+})$ for some $R>0$, that
solves (\ref{eq:gov:tot1}) and (\ref{eq:gov:tot2}), such that 
$u^{\rm tot}-u^{\rm tot}_b$ is radiating. 
Since $u^{\rm tot}-u^{\rm tot}_b$ vanishes on $\Gamma_R^{\rm ext}$, 
the ITBC condition (\ref{eq:tbc:1}) and (\ref{eq:tbc:2}) can still be applied to
obtain a boundary value problem for $u=u^{\rm tot}$ in $\Omega_R^-$, in the
distributional sense,
\begin{align}
  \label{eq:u1:2}
  \Delta u +  k^2 u &= 0,\quad{\rm on}\quad\Omega_R^-,\\
  \label{eq:u2:2}
  \gamma u &= 0,\quad{\rm on}\quad\Gamma_{\rm int},\\
  \label{eq:u3:2}
  2\gamma u &= \cS \phi + \bi \eta \cK(E^*\cS_0^2 E\phi) + \bi \eta E^*\cS_0^2E\phi + 2\gamma u^{\rm tot}_b,\\
  \label{eq:u4:2}
  2\partial_{\nu} u &= \cKd \phi - \phi + \bi \eta \cT(E^*\cS_0^2 E\phi) + 2\partial_{\nu} u^{\rm tot}_b.
\end{align}

Similar to Problem (P1), we can pose an equivalent variational formulation as
follows: Find $(u,\phi)\in V$, such that, for any $(v,\psi)\in V$, 
\begin{equation}
  \noindent{\rm (P2):}\quad\quad a((u,\phi),(v,\psi))=<\partial_{\nu} u^{\rm tot}_b,\gamma v>_{\Gamma_R} + <\gamma u_b^{\rm tot},\psi>_{\Gamma_R}/2.
\end{equation}
We have the following well-posedness result.
\begin{mytheorem}
  \label{thm:plane}
  For any $k>0$ and any incident plane wave $u^{\rm inc}= e^{\bi k(\cos\theta x_1 - \sin\theta
    x_2)}$ with $\theta\in (0,\pi)$, there exists a unique solution $u^{\rm
    tot}\in H^1_{\rm loc}(\Omega)$ such that $u^{\rm tot}-u^{\rm tot}_b$ radiates at infinity.
  \begin{proof}
    The proof follows from that the r.h.s of (P2) defines a bounded functional in $V^*$.
  \end{proof}
\end{mytheorem}
In practice, it is extremely expensive to evaluate $u^{\rm tot}_b$.
Nevertheless, we have pointed out in \cite{luhu19} that by merely extracting the
plane waves $u^{\rm inc}+u_{\pm}^{\rm ref}$ in $\Omega \cap
\Omega_{b,\theta}^{\pm}$ from $u^{\rm tot}$, respectively, we get 
\begin{equation}
  \label{eq:def:u_d}
  u_d = \left\{
    \begin{array}{lc}
      u^{\rm tot}-u^{\rm inc}- u^{\rm ref}_{-},\quad \Omega\cap \Omega_{b,\theta}^-,\\
      u^{\rm tot} - u^{\rm inc} - u^{\rm ref}_{+},\quad \Omega\cap \Omega_{b,\theta}^+,
    \end{array}
  \right.
\end{equation}
piecewisely satisfying an integral form of SRC condition. In fact, the following
lemma indicates that $u_d$, though discontinuous across $L_\theta$, satisfies
the stronger uSRC condition (\ref{eq:src:1}) and (\ref{eq:src:2}).
\begin{mytheorem}
  \label{thm:plane:2}
  For any incident wave $u^{\rm inc}(x)=e^{\bi k(\cos\theta x_1 - \sin\theta
    x_2)}$ with incident angle $\theta\in(0,\pi)$, $u_d$ defined in
  (\ref{eq:def:u_d}) radiates in $\Omega\cap\Omega_{b,\theta}^\pm$ as
  follows: if $x_2\geq 0$, so that
  $x=(r\cos\alpha,r\sin\alpha)$ with $\alpha\in[0,\pi]$, then we have,
  \begin{equation}
    \label{eq:asym:ud}
    \partial_{r} u_d({x}) - \bi k u_d(x) = {\cal O}(r^{-3/2}),\quad{\rm as}\quad r\to\infty,
  \end{equation}
uniformly for $\alpha\in[0,\pi]$; if $x=(x_1,x_2)$ for $x_2\in [-h,0]$, then
  we have
  \begin{equation}
    \label{eq:asym:ud2}
    \partial_{x_1} u_d({x}) - \bi k u_d(x) = {\cal O}(x_1^{-3/2}),\quad{\rm
      as}\quad x_1\to+\infty,
  \end{equation}
  uniformly for $x_2\in[-h,0]$. 
  \begin{proof}
    Suppose $x=(r\cos\alpha,\sin\alpha)\in\mathbb{R}_+^2$ for $\alpha\in
    [0,\pi]$ first. Let $\hat{x}=(\cos\alpha,\sin\alpha)$ and
    $\hat{x'}=(\cos\theta,\sin\theta)$ so that $\nu'=(\sin\theta,-\cos\theta)$.
    By Lemma~\ref{lem:backsol} and Theorem~\ref{thm:plane}, it is sufficient to
    prove that
    \[
      u^{\rm tot}_b(x) - u^{\rm inc}(x) - u^{\rm ref}_{\pm}(x)=(1-2
      e^{2\bi \sin\theta h})\int_{L_\theta} \partial_{\nu(x')}G(x;x')e^{\bi k
        x'\cdot(\cos\theta,\sin\theta)}ds(x'),
    \]
    satisfies (\ref{eq:asym:ud}) and (\ref{eq:asym:ud2}), as $u^{\rm tot}_b -
    u^{\rm tot}$ is radiating. As $x'=r'\hat{x'}\in L_\theta$, we can decompose
    the above integral as
    \begin{align}
      \label{eq:def:I123}
      \int_{0}^{\infty} \partial_{\nu(x')}G(x;x') e^{\bi k r'} dr'=&\int_{0}^{\infty} \partial_{\nu(x')}G_1(x;x') e^{\bi k r'} dr' + \int_{0}^{\infty} \partial_{\nu(x')}\Phi_k(x;x') e^{\bi k r'} dr' \nonumber\\
      &+ \int_{0}^{+\infty} \partial_{\nu(x')}(-\Phi_k(x;x_{\rm im}')) e^{\bi k r'} dr':=I_1(x) + I_2(x) + I_3(x),
    \end{align}
    where $x_{\rm im}'=(x',-y')$ is the imaging point of $x'=(x',y')$ about $y=0$. 
    We consider $I_1$ first. By (\ref{eq:sol:G1}), we get
    \[
      I_1(x) = \frac{1}{2\pi}\int_{0}^{+\infty}\int_{\cal L}\partial_{\nu(x')}\hf^+(\xi;x')e^{-\bi \xi r\cos\alpha +\bi\mu r\sin\alpha} d\xi e^{\bi k r'} d r'.
    \]
    As indicated in the proof of Lemma~\ref{lem:G1:case1}, we can deform ${\cal
      L}$ to a proper path $\tilde{\cal L}$ so that $\tilde{\cal L}\cap {\cal
      L}=\emptyset$ and that $\tilde{\cal L}$ contains a small neighborhood of
    the single stationary point $\xi=-k\cos\alpha$ on the real axis.
    Then, one could use the same asymptotic analysis as in the proof of
    Lemma~\ref{lem:G1:case1} for the large argument $r'$ of $\hf^+$ in
    (\ref{eq:+t-}) to prove, for any $\xi\in \tilde{\cal L}$,
    $\partial_{\nu(x')}\hf^+(\xi;x')$ is smooth of $\xi$ and decays of order
    ${\cal O}(r'^{-3/2})$ as $r'\to\infty$. Exchanging the order of integration
    and using the same asymptotic analysis again but for the large argument $r$
    in $I_1$, we get $\partial_r I_1({x}) - \bi k I_1(x) = {\cal O}(r^{-3/2})$,
    as $r\to\infty$, uniformly for $\alpha\in[0,\pi]$.

    The most troublesome term is $I_2$ with logarithmically singular kernel.
    Observing that
    \[
      \partial_r = \partial_{\hat{x}} = [(\hat{x}\cdot
      \hat{x'})\hat{x'} + (\hat{x}\cdot\nu')\nu']\cdot\nabla_{x},
    \]
    and $\nabla_{x} \Phi_k(x;x') = -\nabla_{x'}\Phi_k(x;x')$, we get
    \begin{align*}
      \left(\partial_r- \bi k\right) I_2(x) =&- \int_{0}^{+\infty} \left[ \partial_{\nu'}\partial_{ r' }\Phi_k(x;x')(\hat{x}\cdot x')+ \partial^2_{\nu'^2}\Phi_k(x;x')(\hat{x}\cdot\nu') + \bi k \partial_{\nu'}\Phi_k(x;x') \right]e^{\bi k r'}d r'.
    \end{align*}
    Since for $x \neq x'$, $\partial^2_{\nu'^2}\Phi_k = -\partial^2_{r'^2}\Phi_k - k^2\Phi_k$, integration by parts gives
    \begin{align*}
      \partial_r I_2({x}) - \bi k I_2(x) =\partial_{\tau(x)} \Phi_k(x;{\bf 0})+ \bi k \Phi_k(x;{\bf 0})(\hat{x}\cdot\nu') + \bi k (\hat{x}\cdot \hat{x'}-1)I_2(x),
    \end{align*}
    where $\tau(x)= (\sin(\alpha),-\cos(\alpha))$. Thus,
    \begin{align*}
      I_2(x)
      =&\frac{\bi k}{4}\int_{0}^{\infty}e^{-\bi k|x-x'|}H_1^{(1)}(k|x-x'|)\frac{r<\hat{x},\nu'>}{|x-x'|} e^{\bi k (r' + |x-x'|)}dr'\\
      =& \frac{\bi k}{4}\int_{1}^{\infty}e^{-\bi k rf(t)}H_1^{(1)}(k rf(t))\frac{r<\hat{x},\nu'>}{t-<\hat{x},\hat{x'}>} e^{\bi k rt}dt,
    \end{align*}
    where we have introduced a new variable $t=(r' + |x - x'|)/r$ with $t\geq 1$
    and have used the following identities.
    \begin{align*}
      r'&=\frac{rt^2-r}{2t-2\cos(\theta-\alpha)},\quad f(t)=\frac{t^2-2t\cos(\theta-\alpha)+1}{2t-2\cos(\theta-\alpha)},\\
      |x-x'| &= rf(t),\quad \frac{dr'}{dt} = \frac{rf(t)}{t-\cos(\theta-\alpha)}>0.
    \end{align*}
    For any $x\neq 0$ \cite{nist10},
    \[
      \left|e^{-\bi x}H_1^{(1)}(x) - e^{-\bi 3\pi/4}\sqrt{\frac{2}{\pi x}}(1 +
        \frac{3\bi}{8x})\right|\leq C x^{-5/2},
    \]
    for some constant $C>0$. Integration by parts,
    \begin{align*}
      I_2^1(x):=&\frac{\bi k}{4}\int_{1}^{\infty}e^{-\bi 3\pi/4}\sqrt{\frac{2}{\pi krf(t)}}(1 + \frac{3\bi}{8 krf(t)})\frac{r\sin(\theta-\alpha)}{t-\cos(\theta-\alpha)} e^{\bi k rt}dt\\
      =&-\frac{e^{-\bi 3\pi/4}}{4}\sqrt{\frac{2}{\pi kr}}(1 + \frac{3\bi}{8 kr})\frac{\sin(\theta-\alpha)}{1-\cos(\theta-\alpha)} e^{\bi k r}\nonumber\\
                &-\frac{\sin(\theta-\alpha)}{4}e^{-\bi 3\pi/4}\sqrt{\frac{2}{\pi kr}}\int_{1}^{\infty}\left(  f(t)^{-1/2}+ \frac{3\bi f(t)^{-3/2}}{8 kr(t-\cos(\theta-\alpha))}\right)' e^{\bi k rt}dt.
    \end{align*}
    By routine calculations, 
    \begin{align*}
      \int_{1}^{\infty} (f(t)^{-1/2})'e^{\bi k r t}dt &= \frac{e^{\bi k r\cos(\theta-\alpha)}}{\sqrt{2}|\sin(\theta-\alpha)^{1/2}|}\int_{\frac{1-\cos(\theta-\alpha)}{|\sin(\theta-\alpha|)}}^{+\infty}\frac{(t^2-1)e^{\bi k r t|\sin(\theta-\alpha)|}}{t^{1/2}(t^2+1)^{3/2}}dt.
    \end{align*}
    Integration by parts,
    \begin{align*}
      \left|  \int_{1}^{\infty} (f(t)^{-1/2})'e^{\bi k r t}dt\right| \leq&  C |\sin(\theta-\alpha)|^{-2}r^{-1} \\
      &+ |\sin(\theta-\alpha)|^{-3/2}r^{-1}\int_{\frac{1-\cos(\theta-\alpha)}{|\sin(\theta-\alpha|)}}^{+\infty}\left| \left( \frac{t^2-1}{t^{1/2}(t^2+1)^{3/2}} \right)' \right|dt\\
      \leq& C |\sin(\theta-\alpha)|^{-2}r^{-1},
    \end{align*}
    for constant $C$ independent of $\alpha$ and $r$.
    Next,
    \begin{align*}
 &\left|  \int_{1}^{\infty} \left(\frac{f(t)^{-3/2}}{t-\cos(\theta-\alpha)}\right)'e^{\bi k r t}dt\right| \\
=& \left|  \frac{2e^{\bi k r\cos(\theta-\alpha)}}{|\sin(\theta-\alpha)|^{5/2}}\int_{\frac{1-\cos(\theta-\alpha)}{|\sin(\theta-\alpha|)}}^{+\infty}\left[ \frac{t^{-1/2}}{2(t^2+1)^{3/2}} - \frac{4t^{3/2}}{3(t^2+1)^{5/2}}\right]e^{\bi k r |\sin(\theta-\alpha)| t}dt\right|\\
\leq &C |\sin(\theta-\alpha)|^{-5/2}
    \end{align*}
    so that one gets
    \begin{align*}
      I_2^1(x) = -\frac{e^{-\bi 3\pi/4}}{4}\sqrt{\frac{2}{\pi k r}}\frac{\sin(\theta-\alpha)}{1-\cos(\theta-\alpha)}e^{\bi k r} + {\cal O}(r^{-3/2}\sin^{-3/2}(\theta-\alpha)).
    \end{align*}
    On the other hand, 
    \begin{align*}
      |I_2(x) - I_2^1(x)|&\leq C\int_{1}^{\infty}(k rf(t))^{-5/2}\frac{r|\sin(\theta-\alpha)|}{t-\cos(\theta-\alpha)}dt\\ 
&=Cr^{-3/2}|\sin(\theta-\alpha)|^{-3/2}\int_{\frac{1-\cos(\theta-\alpha)}{|\sin(\theta-\alpha|)}}^{+\infty}\frac{t^{3/2}}{(t^2+1)^{5/2}}dt \leq C r^{-3/2}|\sin(\theta-\alpha)|^{-3/2}.
    \end{align*}
    Consequently, from the fact that 
    \[
      \partial_{\tau(x)}\Phi_k(x;{\bf 0}) = {\cal O}(r^{-3/2}),\quad
      \Phi_k(x;{\bf 0}) = \frac{\bi}{4}\sqrt{\frac{2}{\pi k
          r}}e^{\bi(kr-\pi/4)} + {\cal O}(r^{-3/2}),
    \]
    and from the estimates for $I_2^1(x)$ and $I_2(x)-I_2^1(x)$ above, we get
    $\partial_r I_2 - \bi k I_2 = {\cal O}(r^{-3/2})$, as $r\to\infty$,
    uniformly for all $\alpha\in [0,\pi]$. As for $I_3$, since the kernel
    function is smooth everywhere, routine calculations, analogous to but much
    simpler than that of $I_2$, give $\partial_r I_3 - \bi k I_3 = {\cal
      O}(r^{-3/2})$, as $r\to\infty$, uniformly for all $\alpha\in [0,\pi]$,
    which concludes the proof for $y>0$. The case when $y\in [-h,0]$ is much
    easier to prove; we omit the details.
  \end{proof}
\end{mytheorem}
\begin{myremark}
  \label{rm:ud}
  From the proof, one obtains the following asymptotic behavior of
  $I_2$ in (\ref{eq:def:I123}): For $x=(r\cos\alpha,r\sin\alpha)$ with
  $\alpha\in[0,\pi]$
   \[
     I_2(x) = \frac{e^{\bi k r}}{\sqrt{r}}\left[ \frac{-e^{-\bi
           3\pi/4}\sin(\theta-\alpha)}{4(1-\cos(\theta-\alpha))}\sqrt{\frac{2}{\pi
           k}} + {\cal
         O}\left(\frac{\sin(\alpha-\theta)^{-3/2}}{r}\right)\right],\ {\rm
     as}\ r\to\infty,
   \]
   so that the far-field pattern of $I_2(x)$ is finite only when the observation
   angle $\alpha$ is away from the reflective angle $\theta$, and approaches
   $\infty$ as $\alpha\to\theta$; roughly speaking, $I_2$ and hence $u_d(x)$,
   though still satisfying the uSRC condition (\ref{eq:src:1}) and
   (\ref{eq:src:2}), don't have a far-field pattern uniformly for
   $\alpha\in[0,\pi]$ unless a small neighborhood of $\theta$ is removed.
 \end{myremark}

\section{Numerical experiments}
As the background Green function $G$ is quite expensive to evaluate, the ITBC
developed in this paper will not be used to numerically truncate the unbounded
domain $\Omega$. Instead, we shall truncate $\Omega$ by the well-known perfectly
matched layer (PML) \cite{ber94,chewee94}.

\subsection{The PML-BIE method}
Mathematically, PML corresponds to complexified transformations of the axial
variables, i.e.,
\begin{equation}
  \label{eq:pml:trans}
  \tilde{x}_j = x_j + \bi \int_0^{x_j} \sigma_j(t)dt,
\end{equation}
where the absorbing function $\sigma_j(t)\geq 0$ is strictly positive only in
the interval $|x_j|>L_j/2$, known as the PML layer, for $j=1,2$. In
\cite{luhu19}, a numerical mode matching method was developed for
$\partial\Omega$ restricted to involving rectangular interfaces only. Here, we
shall adopt the high-accuracy PML-BIE method developed in \cite{luluqia18},
since it allows more general scattering surfaces. In the following we shall
consider plane incidences only; the case for cylindrical incidences can be
formulated similarly but more easily.

As illustrated before, the unknown radiating wave $u=u^{\rm tot}-u^{\rm tot}_b$
seems to be more suitable to compute in practice since it decays exponentially
in the PML. However, it requires calculating the time-consuming background
solution $u_b^{\rm tot}$. Following our previous work \cite{luhu19}, we choose
to compute the discontinuous and radiating wave $u_d$ in \eqref{eq:def:u_d}
(c.f. Remark \ref{rm:ud}) as only simple plane waves are subtracted from $u^{\rm
  tot}$. However, a pseudointerface, i.e., ray $L_\theta$, arises naturally
due to the following jumping conditions across $L_\theta$
\begin{align}
 \label{eq:cond:ud:1} 
 u_d^+(x) - u_d^-(x) &= u^{-}_{\rm ref}(x) - u^+_{\rm ref}(x),\\ 
 \label{eq:cond:ud:2} 
 \partial_{\rm \nu} u_d^+(x) &= \partial_{\rm \nu} u_d^-(x),
\end{align}
where $\nu=(\sin\theta,-\cos\theta)$ is the unit normal vector to $L_\theta$
pointing rightward, and the superscript $\pm$ indicates the sided-limit taken
from $\pm\infty$ on the $x_1$-axis. 

\begin{figure}[!ht]
  \centering
  \includegraphics[width=0.3\textwidth]{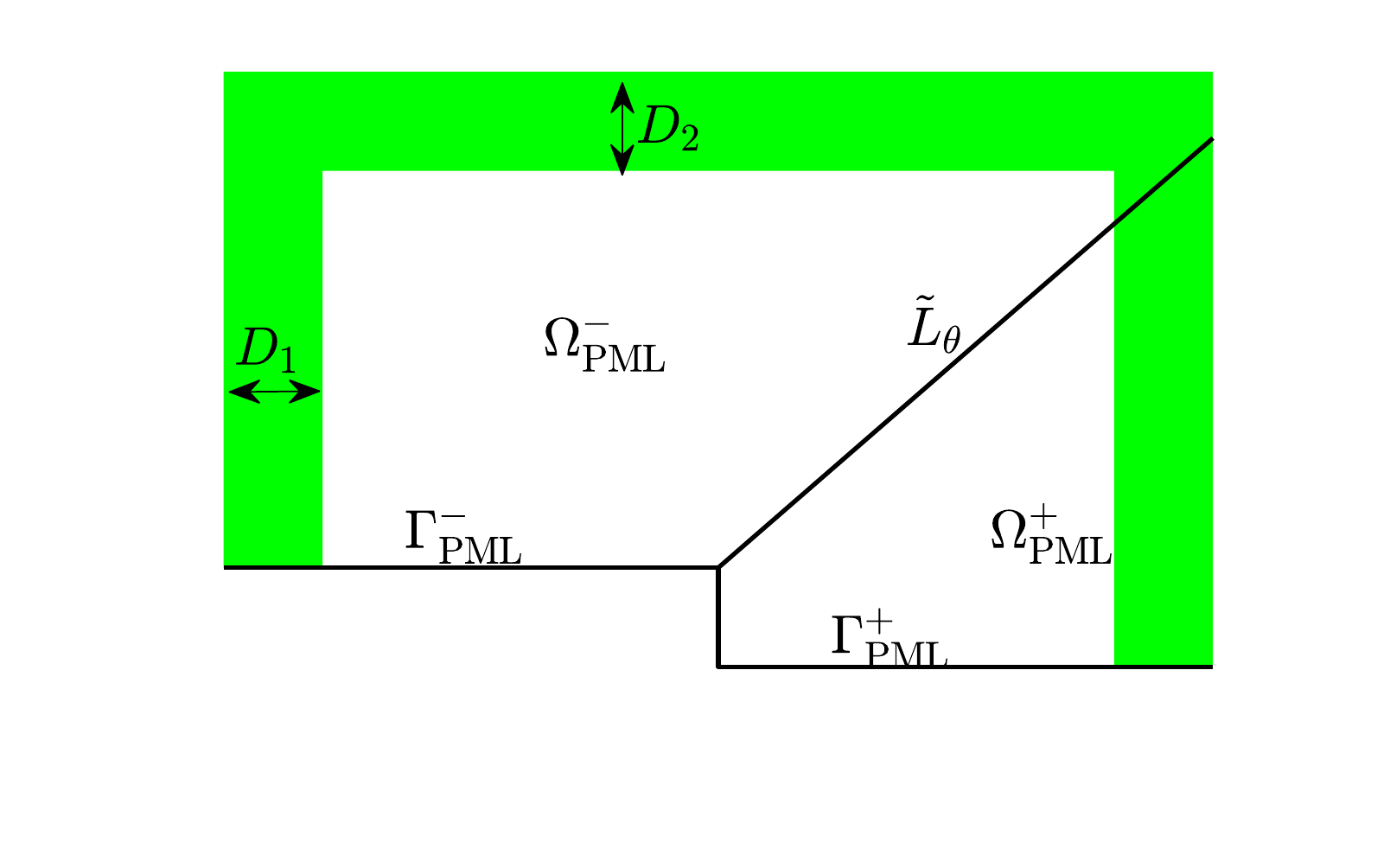}
  \caption{Configuration of the PML. The segment $\tilde{L}_\theta$ is parallel
    to reflective direction for incident angle $\theta$. The truncated domain
    $\Omega_{\rm PML}$ is split into $\Omega_{\rm PML}^{\pm}$ by
    $\tilde{L}_\theta$, while the truncated surface $\Gamma_{\rm PML}$ is split into $\Gamma_{\rm
      PML}^{\pm}$.}
  \label{fig:modelpml}
\end{figure}

We summarize the basic idea of the PML-BIE method \cite{luluqia18} here. The
associated PML configuration for scattering surface $\partial\Omega_h$ is shown in
Figure~\ref{fig:modelpml}.
Let $\Omega_{\rm PML} = \Omega\cap
(-L_1/2-D_1,L_1/2+D_1)\times(-L_2/2-D_2,L_2/2+D_2)$ be the truncated domain for
PML thickness $D_j$ in $x_j$-direction, ${\Omega}^{\pm}_{\rm PML}$ be the two
sub-domains of $\Omega_{\rm PML}$ separated by $\tilde{L}_\theta=L_\theta\cap
\Omega_{\rm PML}$, $\Gamma_{\rm PML} = \partial\Omega\cap \overline{\Omega_{\rm
    PML}}$ be the truncated surface and $\Gamma^{\pm}_{\rm PML}$ be the two
sub-surfaces of $\Gamma_{\rm PML}$ separated by $\tilde{L}_\theta$. Then, with
the coordinate transformations (\ref{eq:pml:trans}), $\tilde{u}_d(x_1,x_2):=
u_d(\tilde{x}_1,\tilde{x}_2)$ satisfies
\begin{align}
  \label{eq:gov:tud}
\nabla\cdot({\bf A}\nabla \tu_d) + k^2J \tu_d &= 0,\quad{\rm on}\quad\Omega_{\rm PML}^{\pm},\\
  \label{eq:bc:tud}
\tu_d&= 0,\quad{\rm on}\quad\partial \Omega_{\rm PML}\backslash\overline{\Gamma_{\rm PML}},
\end{align}
with the following interface conditions on $\Gamma_{\rm PML}^{\pm}$ and across
$\tilde{L}_\theta$,
\begin{align}
  \label{eq:cond:tud:pml:1} 
\tu_d(x_1,x_2)&= u^{\rm inc}(\tx_1,\tx_2) + u_{\pm}^{\rm ref}(\tx_1,\tx_2),\quad{\rm on}\quad\Gamma_{\rm PML}^{\pm},\\
  \label{eq:cond:tud:pml:2} 
 \tilde{u}_d^+(x_1,x_2) - \tilde{u}_d^-(x_1,x_2) &= u^{-}_{\rm ref}(\tx_1,\tx_2) - u^+_{\rm ref}(\tx_1,\tx_2),\quad{\rm on}\quad \tilde{L}_\theta,\\
  \label{eq:cond:tud:pml:3} 
\partial_{{ \nu }^c} u_d^+(x_1,x_2) - \partial_{{ \nu }^c} u_d^-(x_1,x_2) &= \partial_{{ \nu }^c}u^{-}_{\rm ref}(\tx_1,\tx_2) - \partial_{{ \nu }^c}u^+_{\rm ref}(\tx_1,\tx_2),\quad{\rm on}\quad \tilde{L}_\theta.
\end{align}
In the above, ${\bf A}(x_1,x_2)={\rm diag}\{(1+\sigma_2(x_2)\bi)/(1+\sigma_1(x_1)\bi),
(1+\sigma_1(x_1)\bi)/(1+\sigma_2(x_2)\bi)\}$, ${\rm \nu}^c= {\bf A}^{T}\nu$ is
the co-normal vector and $\partial_{{\rm \nu}^c}=\nabla\cdot \nu^c$. Note that,
the continuity of $\partial_{\nu}u_d$ across ${L}_\theta$ in (\ref{eq:cond:ud:2})
breaks down after the PML transformation, i.e., $\partial_{\nu^c}u_{\rm ref}^-\neq \partial_{\nu^c}u_{\rm ref}^+$ in general in (\ref{eq:cond:tud:pml:3}).

As $\tilde{\Phi}^k(x;y)=\Phi^k(\tx;\ty)$ is the fundamental solution of
(\ref{eq:gov:tud}) (c.f. \cite[Thm. 2.8]{lassom01}), we could use $\tilde{\Phi}$ to build up boundary integral
equations for $\tu_d$ and $\partial_{\nu^c} \tu_d$ on the boundary of
$\Omega_{\rm PML}^{\pm}$ (cf. \cite[Eq. (29)]{luluqia18}).  Adopting the high-order
discretization procedure in \cite[Sec. 4]{luluqia18}, the BIEs,
together with the boundary condition (\ref{eq:bc:tud}) on $\partial\Omega_{\rm
  PML}\backslash\overline{\Gamma_{\rm PML}}$ and the interface conditions (\ref{eq:cond:tud:pml:1}-\ref{eq:cond:tud:pml:3}), can be
discretized by a square linear system for the unknowns of
$\partial_{\nu^c}\tu_d$ on grid points of $\Gamma_{\rm
  PML}^{\pm}\cap\tilde{L}_\theta$ and $\tu_d$ on grid points of
$\tilde{L}_\theta$. Finally, the Green's representation theorem \cite[Eq. (27)]{luluqia18} can be used to compute $\tilde{u}_d$ in $\Omega_{\rm
  PML}^{\pm}$ so that $u_d$ and hence $u^{\rm tot}$ become available in the
physical domain of $\Omega_{\rm PML}$, i.e.,
$(-L_1/2,L_1/2)\times(-L_2/2,L_2/2)$. We refer readers to \cite{luluqia18} for
more details of the PML-BIE method.
\subsection{Results}
We here carry out three experiments to show the exponential convergence of
numerical solutions. In all three examples, we choose the following parameters:
$\lambda=1$ so that $k=2\pi$, the incident angle $\theta=\frac{\pi}{3}$,
$L_1=L_2=5$ and $D_1=D_2=D$. We adopt the following PML absorbing function (cf.
\cite[Eq. (20)]{luluqia18})
\begin{equation}
	\label{eq:sigmaj}
	\sigma_j(x_j) = 
	\left\{
		\begin{array}{lc}
			\frac{2Sf_{1}^{8}}{f_1^{8}+f_2^{8}}, & L_j/2\leq x_1\leq L_j/2+D,\\
			S, & x_1>L_j/2+D,\\
			\sigma_j(-x_1), & x_1 \leq -L_j/2,
		\end{array}
	\right.  
\end{equation}
where
\[
	f_1=\left(\frac{1}{2}-\frac{1}{8}\right)\bar{x}_j^3+\frac{\bar{x}_j}{8}+\frac{1}{2},\quad
	f_2 = 1-f_1,\quad \bar{x}_j = \frac{ x_j - (L_j/2+D)}{D},
\]
and $S>0$ determines the magnitude of $\sigma_j$ so that it can be used to
adjust the absorbing strength for PML absorbing an outgoing wave
\cite{chewee94}. We adopt the discretization scheme of six-order of accuracy as
illustrated in \cite{luluqia18}. To ensure that the discretization error is
sufficiently small so that the truncation error due to PML becomes dominant, we
use $500$ grid points on each smooth segment of $\tilde{L}_\theta\cup
\Gamma_{\rm PML}^{\pm}$. To quantify the truncation error, we regard the
numerical solution for $D=2$ and $S=2$ as our reference solution $u_{\rm
  ref}^{\rm tot}$, and compute the following relative error
\[
  E_{\rm rel} = \frac{||{\bm u}^{\rm tot} - {\bm u}_{\rm ref}^{\rm
      tot}||_{\infty}}{||{\bm u}^{\rm tot}_{\rm ref}||_{\infty}},
\]
as $D$ and $S$ vary, where ${\bm u}_{\rm ref}^{\rm tot}$ denotes the vector of
$u_{\rm ref}^{\rm tot}$ at grid points on the physical part of
$\tilde{L}_{\theta}$, etc..

\begin{figure}[!ht]
  \centering
  (a)\includegraphics[width=0.2\textwidth]{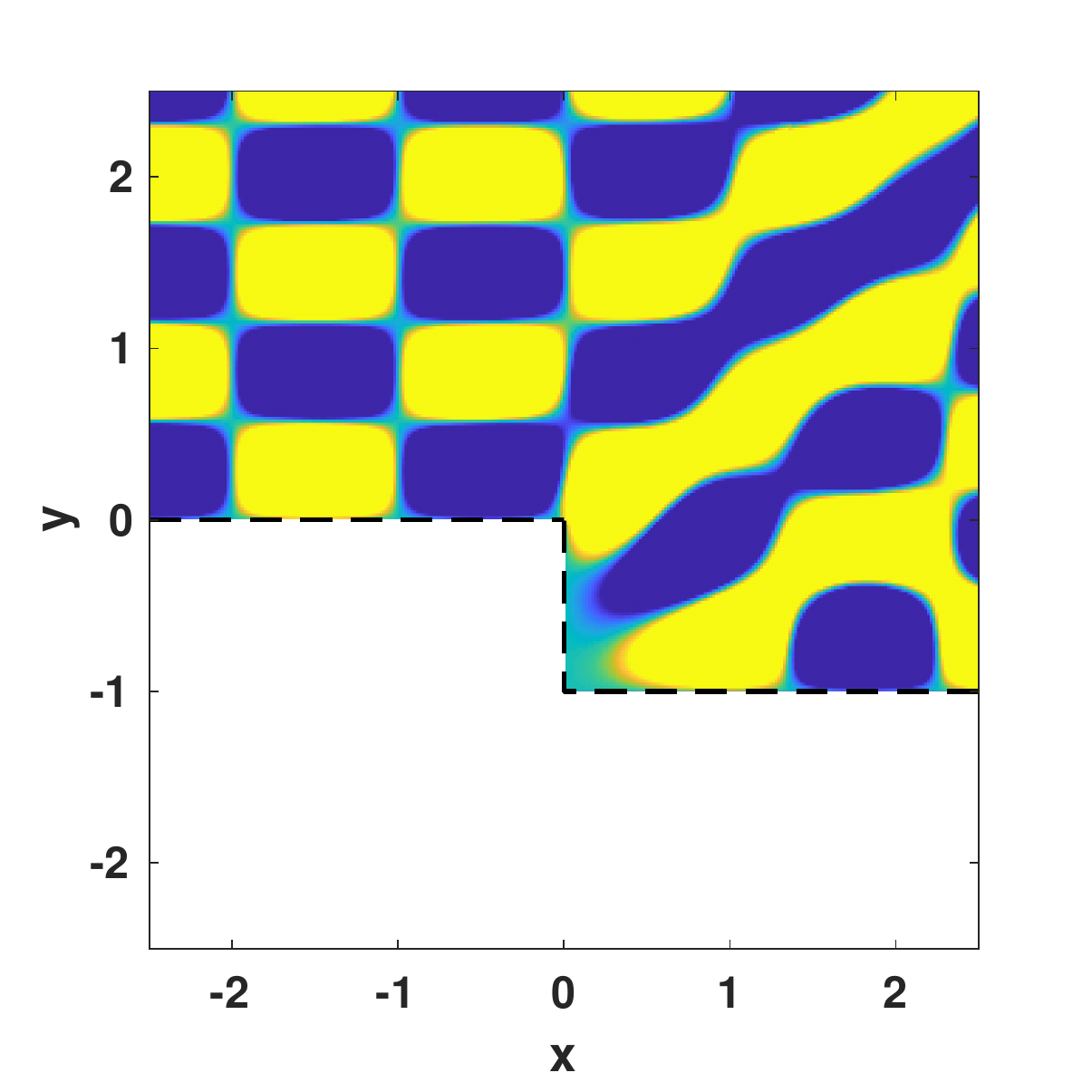}
  (b)\includegraphics[width=0.14\textwidth]{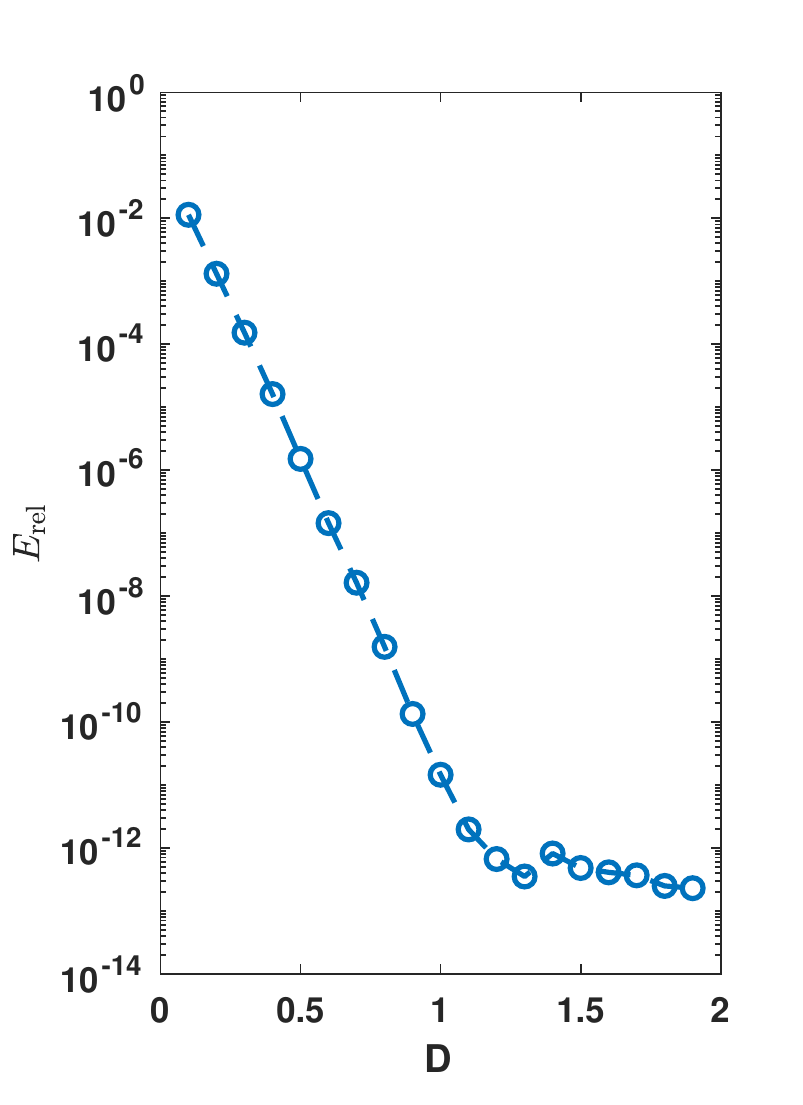}
  (c)\includegraphics[width=0.14\textwidth]{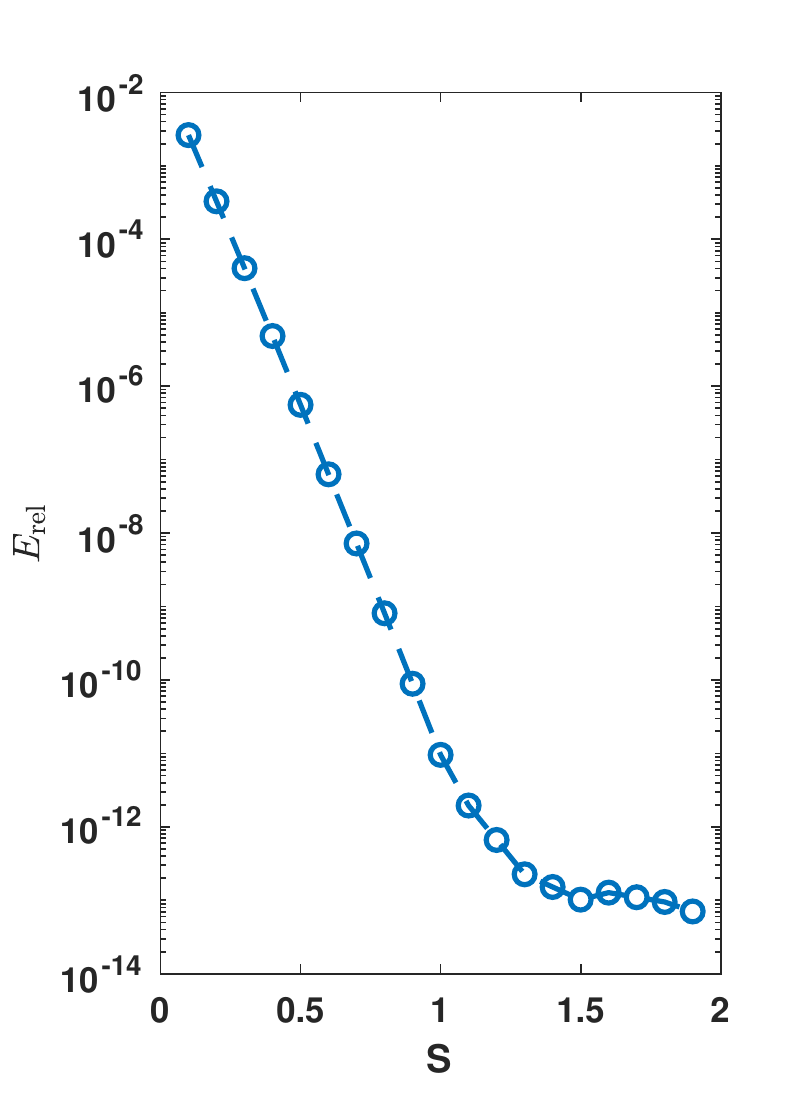}
  \caption{Example 1: (a): real part of the total field $u^{\rm tot}$; (b):
    convergence curve of relative error against PML thickness $D$; (c)
    convergence curve of relative error against PML absorbing constant $S$. The
    dashed lines in (a) indicates the scattering surface $\partial\Omega$.}
  \label{fig:ex1}
\end{figure}

\noindent{\bf Example 1.} The first example has been studied in \cite{luhu19}
where we directly use $\partial\Omega_h$ as the scattering surface. Real part of the reference solution of $u^{\rm tot}$
is plotted in Figure~\ref{fig:ex1}(a),
which is indistinguishable with Figure~4.1(a) in \cite{luhu19}. We show the
convergence curves for $S=2$ and $D$ varying from $0.1$ to $1.9$, and for $D=2$
and $S$ varying from $0.1$ to $1.9$ in Figure~\ref{fig:ex1} (b) and (c),
respectively. We observe that in logarithmic scales, $E_{\rm rel}$ decays
exponentially as $D$ or $S$ increases until  discretization/round-off error
dominates; at least 12 significant digits are obtained by the proposed method.

\begin{figure}[!ht]
  \centering
  (a)\includegraphics[width=0.2\textwidth]{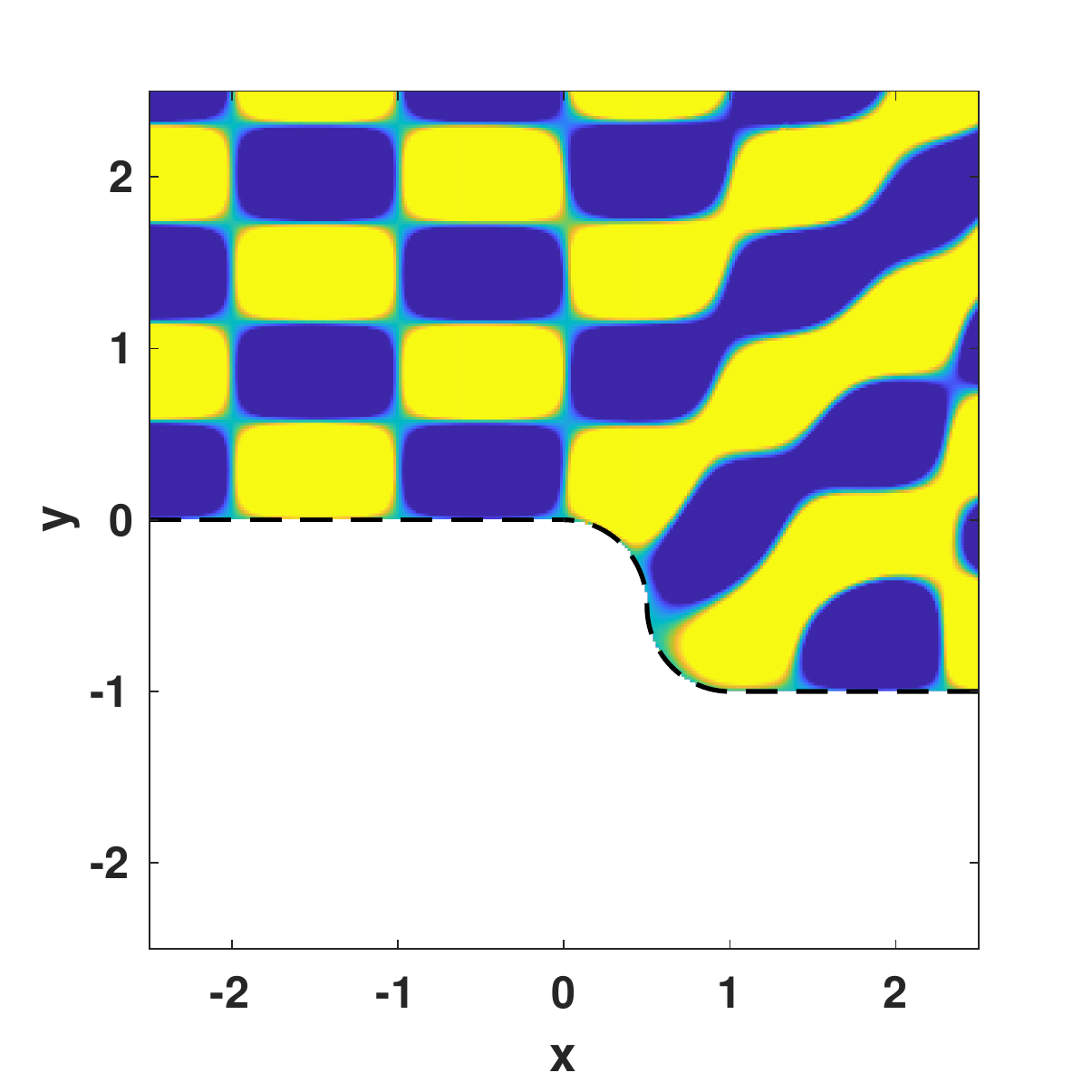}
  (b)\includegraphics[width=0.14\textwidth]{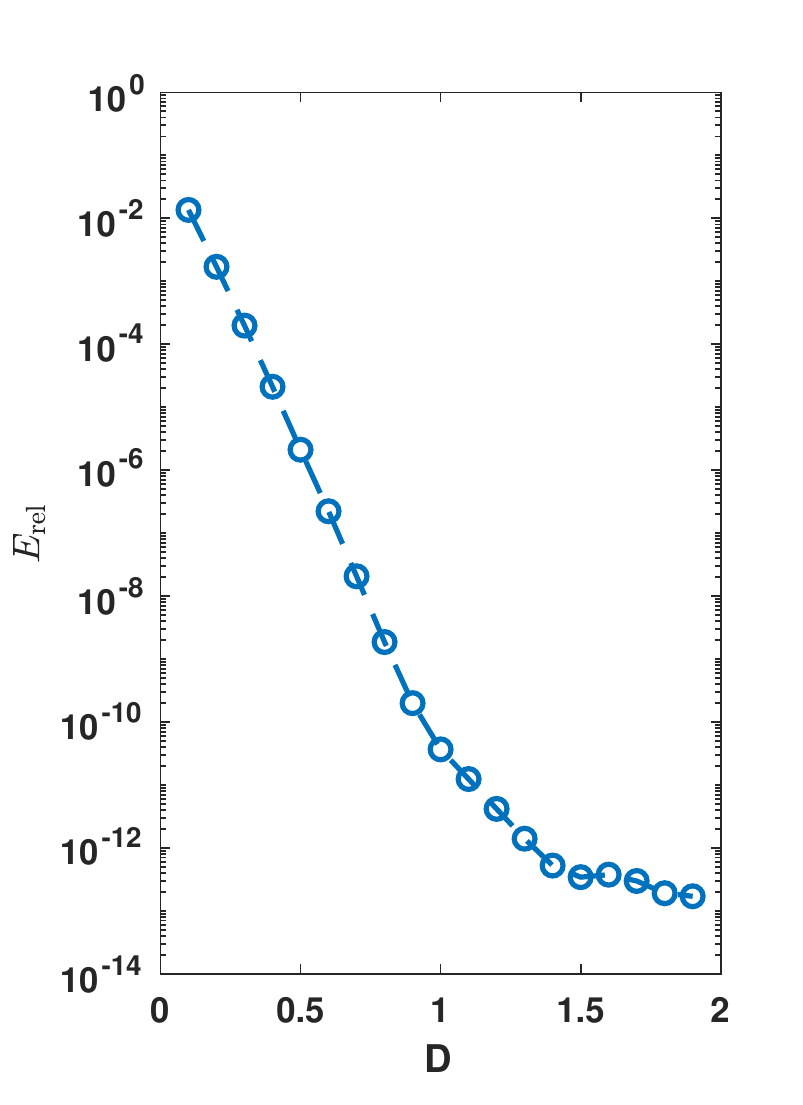}
  (c)\includegraphics[width=0.14\textwidth]{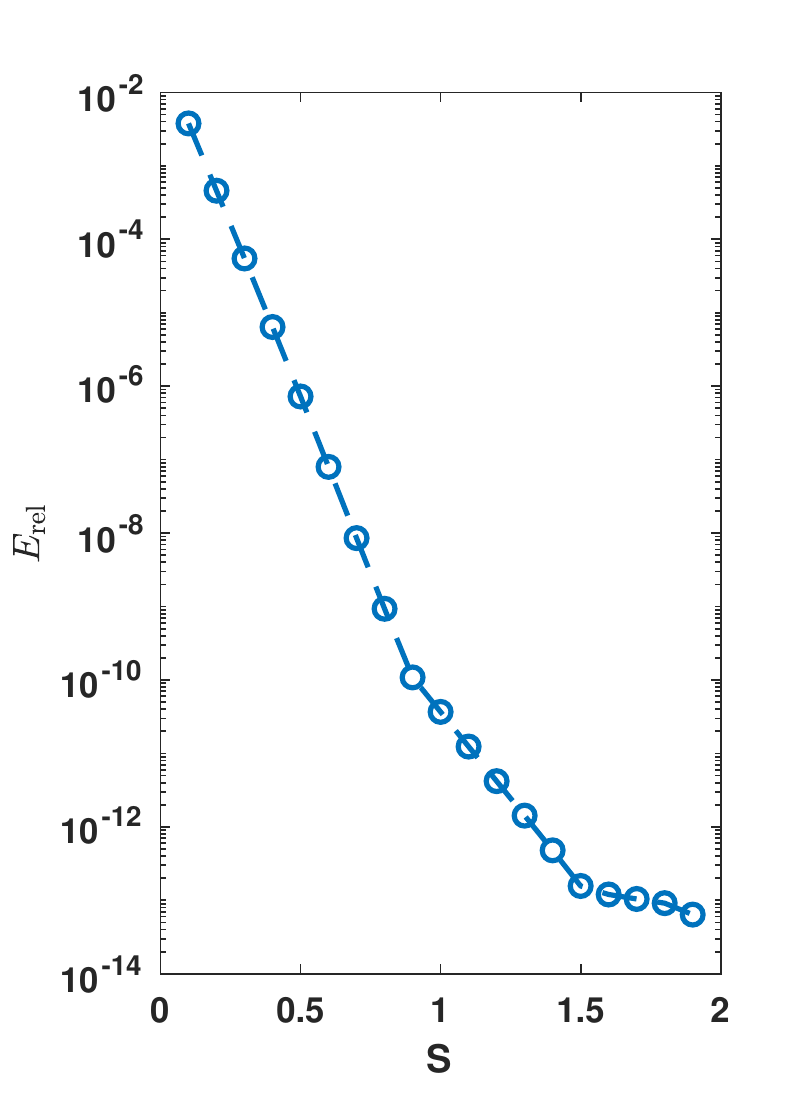}
  \caption{Example 2: (a): real part of the total field $u^{\rm tot}$; (b):
    convergence curve of relative error against PML thickness $D$; (c)
    convergence curve of relative error against PML absorbing constant $S$. The
    dashed lines in (a) indicates the scattering surface $\partial\Omega$.}
  \label{fig:ex2}
\end{figure}

\noindent{\bf Example 2.} In this example, $\Gamma$ consists of two rays and two
quarter circles of radius $1/2$. Real part of the
reference solution of $u^{\rm tot}$ is plotted in Figure~\ref{fig:ex2}(a).
We show the convergence curves for $S=2$ and $D$ varying from $0.1$ to $1.9$,
and for $D=2$ and $S$ varying from $0.1$ to $1.9$ in Figure~\ref{fig:ex2} (b)
and (c), respectively. We observe that in logarithmic scales, $E_{\rm rel}$
decays exponentially as $D$ or $S$ increases until discretization/round-off
error dominates; at least 12 significant digits are obtained by the proposed method. In comparison with Figure~\ref{fig:ex1} (a), $u^{\rm tot}$ in
Figure~\ref{fig:ex2} (a) mainly differs near the two semicircles.

\begin{figure}[!ht]
  \centering
  (a)\includegraphics[width=0.2\textwidth]{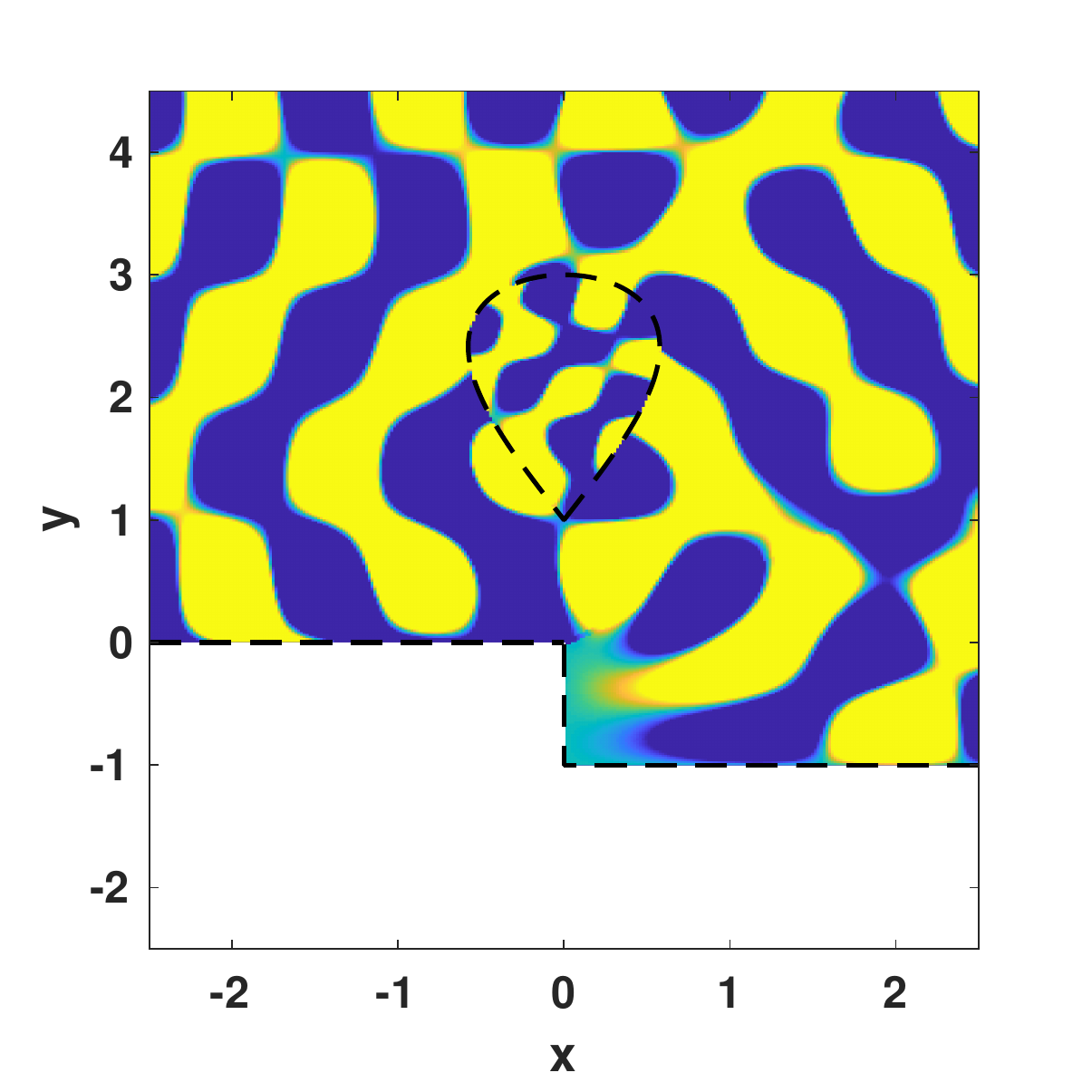}
  (b)\includegraphics[width=0.14\textwidth]{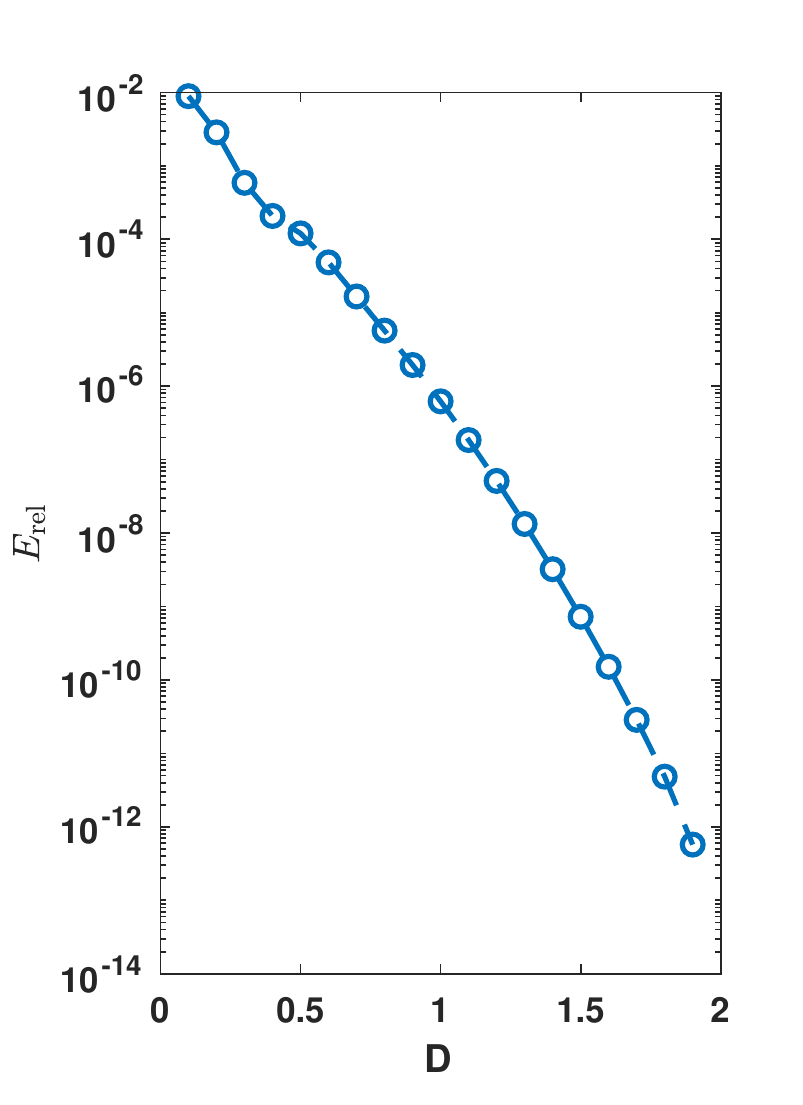}
  (c)\includegraphics[width=0.14\textwidth]{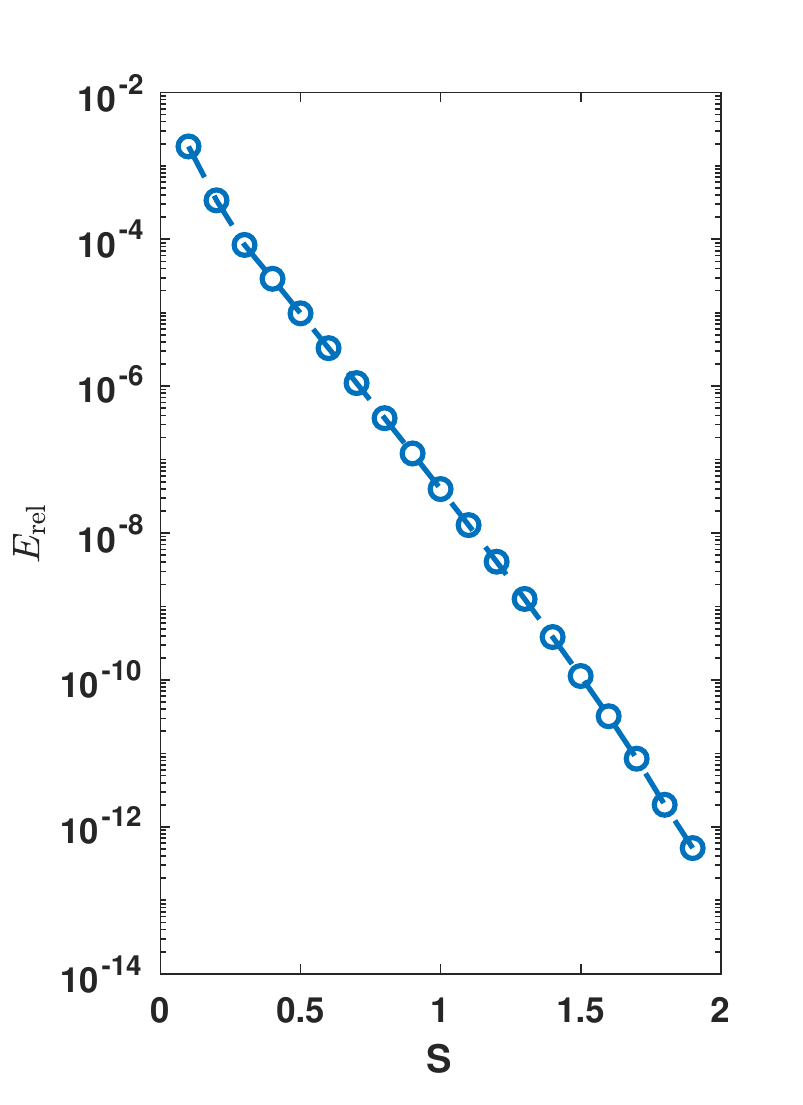}
  \caption{Example 3: (a): real part of the total field $u^{\rm tot}$; (b):
    convergence curve of relative error against PML thickness $D$; (c)
    convergence curve of relative error against PML absorbing constant $S$. The
    dashed lines in (a) indicates the penetrable object and the scattering
    surface $\Gamma$.}
  \label{fig:ex3}
\end{figure}

\noindent{\bf Example 3.} The last example has the same scattering surface
$\partial\Omega_h$ but with a drop-shaped penetrable object of wavenumber
$k_{\rm obj}=2k$ above $\partial\Omega_h$, as
shown in Figure~\ref{fig:ex3}(a).
With a penetrable object, one sees that the corresponding variational
formulation is only a compact perturbation to the problem (P2), so that the
scattering problem still has a unique solution except for a countable set of
wavenumber $k$. Assuming that $k=2\pi$ doesn't lie in this set, we directly use
the PML-BIE method to numerically solve the scattering problem. Real part of the
reference solution of $u^{\rm tot}$ is plotted in Figure~\ref{fig:ex3}(a). We
show the convergence curves for $S=2$ and $D$ varying from $0.1$ to $1.9$, and
for $D=2$ and $S$ varying from $0.1$ to $1.9$ in Figure~\ref{fig:ex3} (b) and
(c), respectively; at least 12 significant digits are obtained by the proposed method.
\subsection{Green function $G$ in the PML and exponential convergence}

\begin{figure}[!ht]
  \centering
  (a)\includegraphics[width=0.3\textwidth]{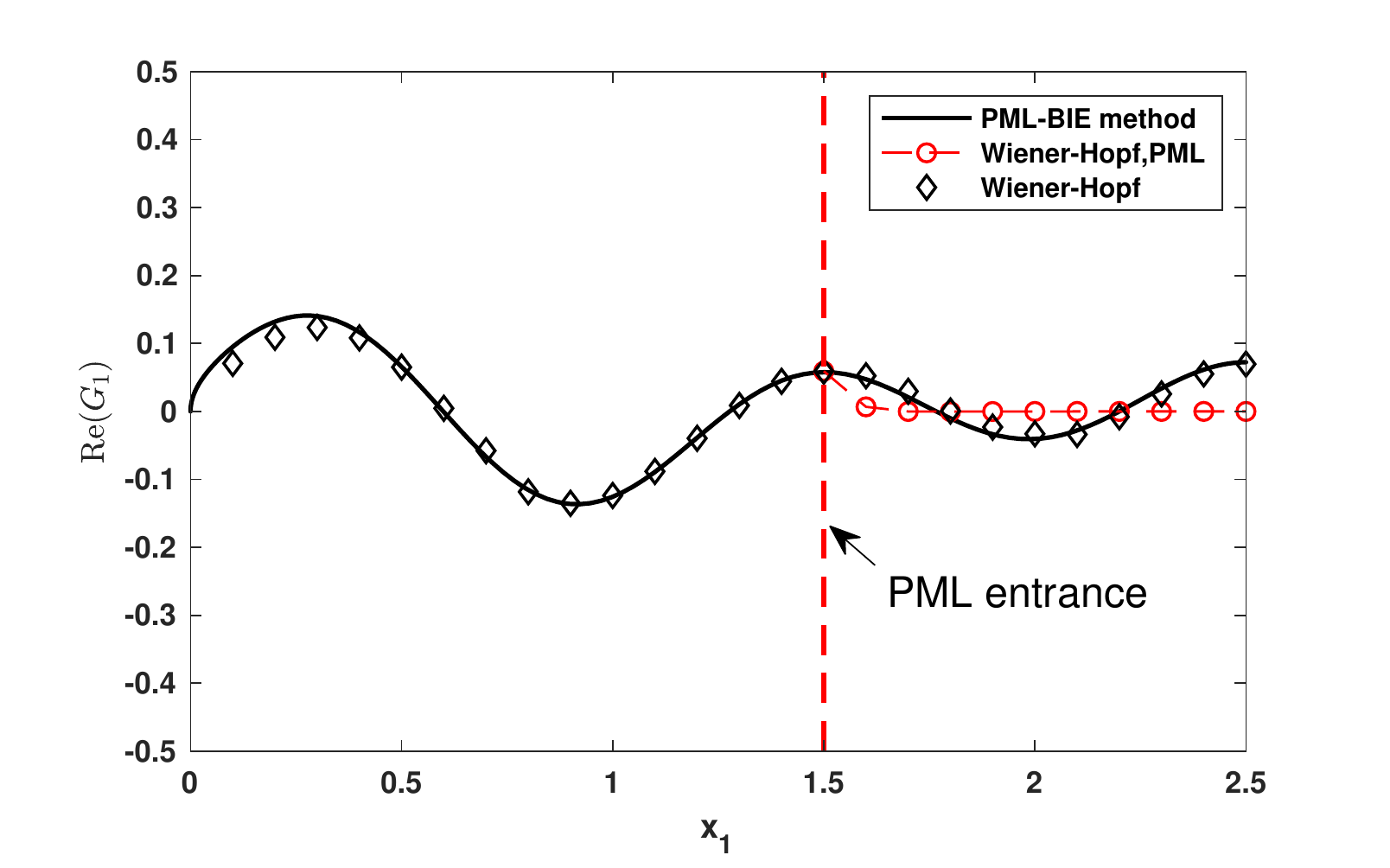}
  (b)\includegraphics[width=0.3\textwidth]{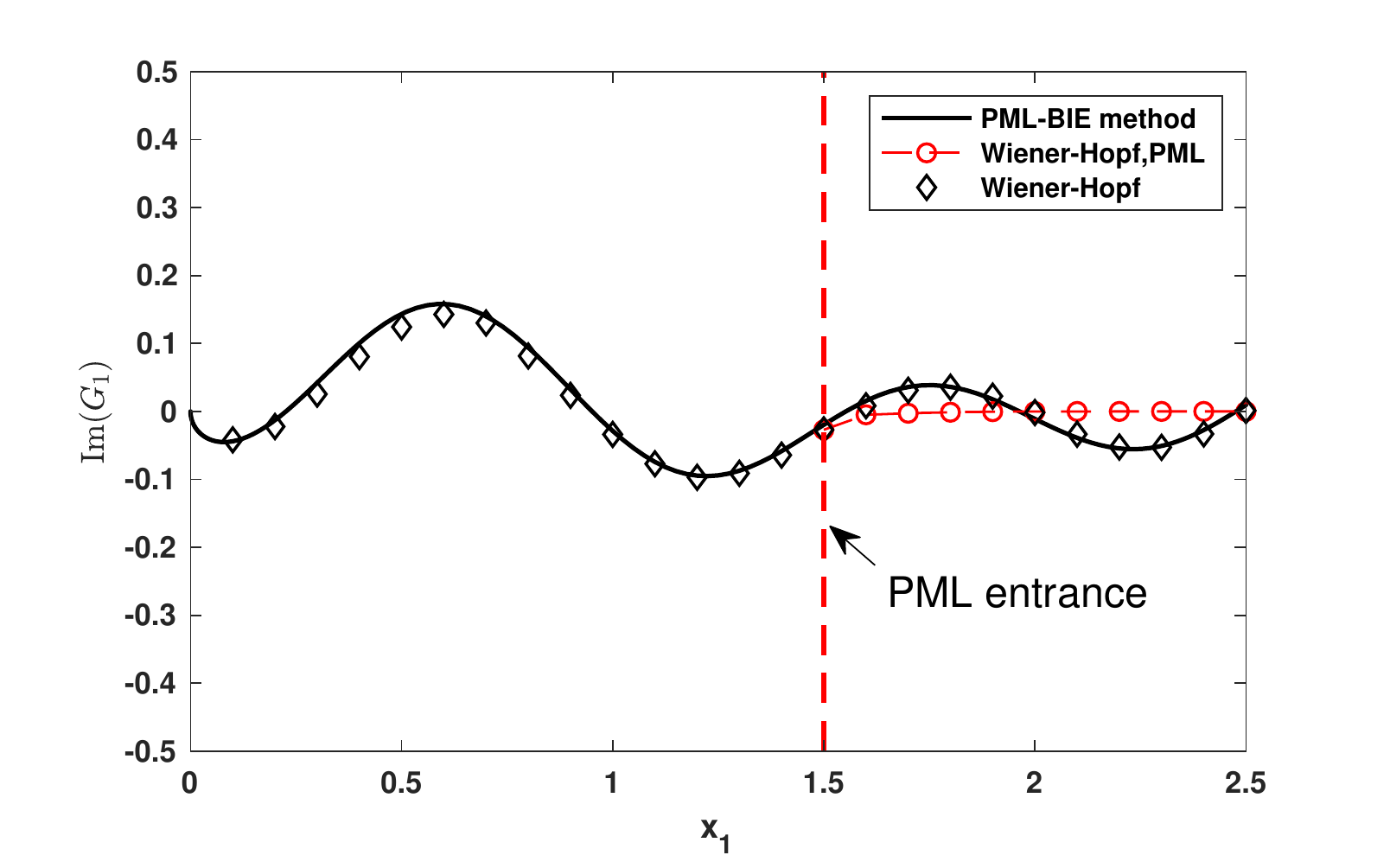}
  \caption{$G_1$ for $k=\frac{2\pi}{1.1}$, $h=1$, source point $x^*=(0,0.4)$ and
    $x=(x_1,0)$ with $x_1\in(0,2.5)$: (a) real part; (b) imaginary part. Solid
    lines represent numerical solution by the PML-BIE method; Diamonds
    ('$\diamond$') represent numerical solution of $G_1(x;x^*)$ by
    (\ref{eq:G1:x1+}); Circle ('$\circ$') dashed lines represent numerical
    solution of $G_1(\tx;x^*)$ in the PML by (\ref{eq:G1+}). The vertical dashed
    lines indicate the entrance of PML.}
  \label{fig:G1}
\end{figure}

We remark that though the background Green function $G$ is of no significant
interest from numerical perspectives, it is the purely outgoing behavior of $G$
in $\Omega$ that makes PML a perfect approach to truncate the unbounded domain
$\Omega$, since the PML truncation error decaying exponentially to 0, as
illustrated by the previous three numerical examples. To conclude this section,
we give a numerical evidence to support this argument and shall defer a rigorous
error analysis in a subsequent paper. Suppose $x^*$ is above the horizontal
axis, when $x=(x_1,x_2)$ is in the PML, we could analytically extend
$G_1(x;x^*)$ by
\begin{equation}
  \label{eq:G1+}
  G_1(\tilde{x};x^*) = \left\{
    \begin{array}{ll}
    \frac{1}{2\pi}\int_{-\infty\to 0\to -\infty\bi}\hf^+(\xi;x^*)e^{-\bi\xi \tx_1+\bi\mu \tx_2}d\xi,& x_1>0,\\
\frac{1}{2\pi}\int_{+\infty \bi\to 0\to +\infty}\hf^+(\xi;x^*)e^{-\bi\xi \tx_1+\bi\mu \tx_2}d\xi, & x_1<0.
    \end{array}
\right.
\end{equation}
where $\tilde{x}=(\tilde{x}_1,\tx_2)$. We could verify numerically that
$G_1(\tx;x^*)$ exponentially decays to $0$ as $x_1, x_2\to\infty$, as
illustrated in Figure~\ref{fig:G1},
where we take $k=\frac{2\pi}{1.1}$, $h=1$, $x^*=(0,0.4)$, $x_2=0$ and
$\tx_1=x_1+3\bi(x_1-1.5)$ for $x_1>1.5$. In Figure~\ref{fig:G1}, we use the
PML-BIE method to compute $G_1((x_1,0);x^*)$ for $x_1\in[0,2.5]$, and directly
use the {\it integral} function in MATLAB 2019a to compute $K^{\pm}$ in
(\ref{eq:def:K+-}), $H^{\pm}$ in (\ref{eq:def:H+-}), and hence $G_1$ in the
physical domain by (\ref{eq:G1:x1+}) and in the PML by (\ref{eq:G1+}). We could
see from Figure~\ref{fig:G1} that numerical solution of PML-BIE method coincides
with that of Wiener-Hopf method to some extent, and that $G_1$ decays
exponentially to $0$ in the PML. A theoretical proof of $G_1$ and $G_2$ decaying
exponentially to $0$ in the PML will be presented in a future work. Thus, the
Green's formula (\ref{eq:greenformula}) implies that any radiating solution $u$
must decay to $0$ exponentially in the PML. Consequently, a PML truncation in
terms of posing zero Dirichlet boundary condition on the boundary of
$[-L_1/2-D_1,L_1/2+D_1]\times[-L_2/2-D_2,L_2/2+D_2]$ with thickness parameters
$D_1$ and $D_2$ is expected to cause a truncation error exponentially decaying
to $0$ as $D_j$ increases, so that the PML-truncated solution must be a
physically correct solution!

\section{Conclusion}
This paper has proposed, for wave propagating in a globally perturbed half plane
with a perfectly conducting step-like surface, the sharper uSRC condition
(\ref{eq:src:1}) and (\ref{eq:src:2}), a closed form of its far-field pattern,
and a high-accuracy numerical solver, based on the Green function of a cracked
half plane. By showing that the Green function asymptotically satisfies the
uSRC, we established a well-posedness theory for either a cylindrical incident
wave due to a line source or a plane incident wave, via an associated
variational formulation through the ITBC (\ref{eq:tbc:1}) and (\ref{eq:tbc:2})
terminating the unbounded domain. The theory reveals that the scattered wave,
post-subtracting a known wave field, inherits the asymptotic behavior of the
background Green function, satisfying the same uSRC. For a plane-wave incidence,
subtracting reflected plane waves, due to non-uniform heights of the step-like
surface at infinity, from the scattered wave in respective regions produces a
discontinuous but outgoing wave satisfying the uSRC. Numerically, we adopted a
previously developed PML-BIE method to solve the problem, demonstrating from
numerical results that the truncation error due to PML decays exponentially fast
as thickness or absorbing power of the PML increases, of which the convergence
relies heavily on the outgoing Green function $G$ decaying exponentially in the
PML.

\appendix
\section{Wiener-Hopf method}
In the appendix, we shall derive for $G(x;x^*)$ when $x_2^*\geq 0$; the case
when $x_2^*\in (-h,0)$ can be analyzed similarly. We consider solving $x_2^*>0$
first; the Green function of source on $\Gamma^+$ is the limit as $x_2^*\to
0^+$, as will be seen in the following.

The scattered field $G_1$ and $G_2$ defined in (\ref{eq:G:1}) solve the
homogeneous Helmholtz equation in its domain of definition and satisfy zero
Dirichlet condition on $\Gamma^-$ and $\Gamma^-\cup \Gamma_h$, respectively.
Furthermore, $G_1$ and $G_2$ satisfy the following interface conditions across
$\Gamma^+$, i.e., for $x_1>0$ and $x_2=0$,
\begin{align}
  \label{eq:bc:int1}
  G_1 - G_2 &= 0,\\
  \label{eq:bc:int2}
  \partial_{x_2} G_1 - \partial_{x_2} G_2 &= g(x_1;x^*) :=-\partial_{x_2}G^{\rm in}((x_1,0);x^*).
\end{align}
Let $f(x_1;x^*)=G_1((x_1,0);x^*)$ for $x_1>0$ with the unknown $f(x_1;x^*)$ to be determined. 
In region $\mathbb{R}_2^+$,  we get from Fourier transforming $G_1$ w.r.t $x_1$, 
\begin{equation}
  \label{eq:hG1}
  \hG_1(x_2;\xi,x^*):= \int_{-\infty}^{+\infty}G_1(x;x^*)e^{\bi \xi x_1} dx_1 = \hat{f}^+(\xi;x^*) e^{\bi \mu x_2},
\end{equation}
where $\hat{f}^+(\xi;x^*) = \int_0^{+\infty}f(x_1;x^*)e^{\bi \xi x_1}dx_1$, so
that
\begin{equation}
  \label{eq:hG1p}
  \hG_1'(0;\xi,x^*)=\bi\mu\hat{f}^+(\xi;x^*).
\end{equation}
Similarly, we get the one-dimensional Fourier transform of $\hG_2$ along $x_1$ 
\begin{equation*}
  \label{eq:hG2}
  \hG_2(x_2;\xi,x^*): = \int_{-\infty}^{+\infty}G_2(x;x^*)e^{\bi \xi x_1} dx_1= \frac{\hat{f}^+(\xi;x^*)}{1-e^{2\bi\mu h}}\left( e^{-\bi \mu x_2}- e^{\bi \mu (x_2+2h)}\right),
\end{equation*}
so that
\begin{equation}
  \label{eq:hG2p}
  \hG_2'(0;\xi,x^*) =-\bi\mu\frac{1+e^{2\bi\mu h}}{1-e^{2\bi \mu h}}\hf^+(\xi;x^*).
\end{equation}
Now denote for $j=1,2$,
\begin{align*}
  \hG_{j,+}(x_2;\xi,x^*) = \int_{0}^{+\infty}G_j(x;x^*)e^{\bi \xi x_1}dx_1,\ {\rm and}\ \hG_{j,-}(x_2;\xi,x^*) = \int_{-\infty}^{0}G_j(x;x^*)e^{\bi \xi x_1}dx_1.
\end{align*}
We get from (\ref{eq:bc:int2}), (\ref{eq:hG1p}) and (\ref{eq:hG2p}) 
the following Wiener-Hopf equation, 
\begin{equation}
  \label{eq:gov:5}
  e^{\bi x_1^*\xi +\bi \mu x_2^*}+[\hG_{1,-}'(0;\xi,x^*)-\hG_{2,-}'(0;\xi,x^*)-g^-(\xi;x^*)] = \frac{2\bi\mu}{1-e^{2\bi\mu h}}\hat{f}^+(\xi;x^*).
\end{equation}
where $\hat{g}^-(\xi;x^*) =\int_{-\infty}^{0}g(x_1;x^*)e^{\bi \xi x_1}dx_1$.
According to (\ref{eq:decomp:1}) and (\ref{eq:decomp:2}), we get 
\begin{equation}
  \label{eq:wh}
  H^-(\xi;x^*)+[\hG_{1,-}'(0;\xi,x^*)-\hG_{2,-}'(0;\xi,x^*)-g^-(\xi;x^*)]\frac{K^{-}(\xi)}{\sqrt{k-\xi}} = \frac{2\bi\sqrt{k+\xi}}{K^+(\xi)}\hat{f}^+(\xi;x^*) - H^+(\xi;x^*).
\end{equation}
In the above, $K^{\pm}$ and $H^{\pm}$ defined in equations (\ref{eq:K+-:onL})
and (\ref{eq:H+-:onL}) can be holomorphically extended to ${\cal C}^{\pm}$, the
two complex regions above and below ${\cal L}$, respectively, via
  \begin{align}
    \label{eq:def:K+-}
    K^{\pm}(\xi) =& \exp\left\{\pm\frac{1}{2\pi\bi}\int_{\cal L} \frac{\Log(1-e^{2\bi\mu(t) h})}{t-\xi}dt\right\},\ {\rm for}\ \xi\in{\cal C}^{\pm}\\
    \label{eq:def:H+-}
    H^{\pm}(\xi;x^*) =& \pm\frac{1}{2\pi\bi}\int_{\cal L} \frac{e^{\bi t x_1^*+\bi \mu(t) x_2^*}K^-(t)}{\sqrt{k-t}(t-\xi)}dt,\ {\rm for}\ \xi\in{\cal C}^{\pm}.
  \end{align}
  According to \cite[p. 33-38]{gak66} and Sokhotski-Plemelj Theorem \cite[Thm.
  7.8]{kre14}, one could obtain the following properties of $K^{\pm}$ and
  $H^{\pm}$:

  (i) $K^{\pm}(\xi)$ is holomorphic in ${\cal C}^{\pm}$ and is
    $C^\infty$-smooth in $\overline{ {\cal C}^{\pm} }$. For sufficiently large
    $|\xi|$ with $\xi\in \overline{ {\cal C}^{\pm} }$, the $m$-th derivative of
    $K^{\pm}(\xi)$ satisfies
    \begin{equation}
      \label{eq:asym:Kpm}
      |(K^{\pm}(\xi))^{(m)} - 1| \leq  C_m(|\xi|^{-\eta}),
    \end{equation}
    for any positive constant $\eta \in(0,1)$ and for any integer $m\geq 0$, where
    $C_m$ is independent of $\xi$.

    (ii) When $x_2^*>0$, $H^{\pm}$ is holomorphic in ${\cal C}^{\pm}$ and is
    $C^\infty$-smooth in $\overline{ {\cal C}^{\pm} }$. For sufficiently large
    $|\xi|$ with $\xi\in \overline{ {\cal C}^{\pm} }$, the $m$-th derivative of
    $H^{\pm}(\xi)$, $|( H^{\pm}(\xi;x^*) )^{(m)}| \leq  C_m|\xi|^{-\eta}$,
    for any positive constant $\eta \in(0,1)$ and for any integer $m\geq 0$, where
    $C_m$ is independent of $\xi$.

    (iii) When $x^*\in\Gamma^+$ so that $x_2^*=0$, Property (ii) for $H^{\pm}$
    doesn't hold since for $t\in{\cal L}$, $\frac{e^{\bi t x_1^*}K^-(t)}{\sqrt{k-t}}$
  is not Holder continuous at infinity so that $H^{\pm}(\xi;x^*)$ as
  $|\xi|=\infty$ may be undefined \cite{gak66}. However, based on the asymptotic
  behavior of $K^-(\xi)-1$ at infinity, one verifies that, by residual theorem,
  defines the following function
\begin{align}
  \label{eq:def:H0pm}
  H^{\pm}(\xi;(x_1^*,0))=& \lim_{x_2^*\to 0^+} H^{\pm}(\xi;x^*) 
  = \pm\frac{1}{2\pi\bi}\left(\int_{\cal L}\frac{e^{\bi t x_1^*}(K^-(t)-1)}{\sqrt{k-t}(t-\xi)}dt + \int_{\cal L}\frac{e^{\bi t x_1^*}}{\sqrt{k-t}(t-\xi)}dt \right)\nonumber\\
      =&\pm\frac{1}{2\pi\bi}\int_{\cal L}\frac{e^{\bi t x_1^*}(K^-(t)-1)}{\sqrt{k-t}(t-\xi)}dt +\left\{
         \begin{array}{lc}
           \frac{e^{\bi \xi x_1^*}}{\sqrt{k-\xi}} & {\rm if}\ \xi\in{\cal C}^+,\\
           0 & {\rm if}\ \xi\in{\cal C}^-.
         \end{array} 
               \right.  
\end{align}
As $\xi$ approach ${\cal L}$ from $ {\cal C}^{\pm}$,
(\ref{eq:H0pm:L}) holds, making $H^{\pm}$ again a $C^\infty$-smooth function in
$\overline{{\cal C}^{\pm}}$, so that for sufficiently large $|\xi|$,
  $| (H^{\pm}(\xi;(x_1^*,0)) )^{(m)}| \leq C_m|\xi|^{-\eta}$,
for any positive constant $\eta \in(0,1/2)$ and for any integer $m\geq 0$, where $C_m$ is independent of $\xi$.

Suppose both sides of (\ref{eq:wh}) vanish at $|\xi|=\infty$ for
$\xi\in\overline{{\cal C}^{\pm}}$. According to Liouville's theorem, we must
have both sides are zero such that (\ref{eq:+t-}) holds, for $\xi\in{\cal L}$.
Finally, from the properties of $K^\pm$ and $H^\pm$ and Morera's theorem, we get properties of
$\hf^+(\xi;x^*)$ in the following lemma.
\begin{mylemma}
  \label{lem:f+} 
  The holomorphic function $\hf^+(\xi;x^*)$ defined in $\overline{{\cal C}^+}$ can be
  extended as a holomorphic function in
  $\mathbb{C}\backslash((-\infty,-k])$ and a continuous function
  in $\overline{\mathbb{C}^{+-}}\cup \overline{\mathbb{C}^{-+}}$. For
  sufficiently large $|\xi|$ with
  $\xi\in\mathbb{R}$, we have $|\hf^+(\xi;x^*)|= {\cal O}(|\xi|^{-1/2-\eta})$,
  where the exponent $\eta$ belongs to $(0,1)$ when $x_2^*>0$ and to
  $(0,\frac{1}{2})$ when $x_2^*=0$.
\end{mylemma}
\section*{Acknowledgment}
The author thanks Profs. Anne-Sophie Bonnet-Bendhia and Christophe
Hazard of ENSTA/POEMS for inspiring this work during the conference of
WAVES 2019, and also thanks Prof. Guanghui Hu of Beijing CSRC, 
Prof. Buyang Li of The Hong Kong Polytechnic University, and Prof. Ya Yan Lu of
City University of Hong Kong for some useful discussions.

\bibliographystyle{plain}
\bibliography{wt}

\end{document}